\documentclass[reqno,12pt]{amsart}
\usepackage{a4wide}

\usepackage{amssymb}

\usepackage{rotating}

\newcommand{\CH}[1][n]{\ensuremath{\mathbb{C}H^{#1}}}
\newcommand{\C}{\ensuremath{\mathbb{C}}}
\newcommand{\R}{\ensuremath{\mathbb{R}}}

\newcommand{\shape}{\ensuremath{\mathcal{S}}}
\newcommand{\pd}{\ensuremath{\mathcal{B}}}
\newcommand{\g}[1]{\ensuremath{\mathfrak{#1}}}
\newcommand{\Id}{\ensuremath{\mbox{Id}}}
\newcommand{\sech}{\mathop{\rm sech}\nolimits}
\newcommand{\csch}{\mathop{\rm csch}\nolimits}
\newcommand{\Exp}{\ensuremath{\mbox{Exp}}}
\newcommand{\enabla}{\ensuremath{\bar{\nabla}}}
\newcommand{\W}[2][\varphi]{\ensuremath{W^{#2}_{#1}}}
\newcommand{\II}{\ensuremath{I\!I}}
\newcommand{\eR}{\ensuremath{\bar{R}}}

\newtheorem{theorem}{Theorem}[section]
\newtheorem*{thm}{Theorem}
\newtheorem{proposition}[theorem]{Proposition}

\newtheorem{corollary}[theorem]{Corollary}

\theoremstyle{remark}
\newtheorem{remark}[theorem]{Remark}

\begin{document}

\title[Homogeneous hypersurfaces in complex hyperbolic spaces]
{\large Homogeneous hypersurfaces\\ in complex hyperbolic spaces}

\author{J\"{u}rgen Berndt, Jos\'{e} Carlos D\'{\i}az-Ramos}

\begin{abstract}
We study the geometry of homogeneous hypersurfaces and their focal
sets in complex hyperbolic spaces. In particular, we provide a
characterization of the focal set in terms of its second
fundamental form and determine the principal curvatures of the
homogeneous hypersurfaces together with their multiplicities.
\end{abstract}


\maketitle

\thispagestyle{empty}

\section{Introduction}

An $s$-representation is the isotropy representation of a
semisimple Riemannian symmetric space. A result by Hsiang and
Lawson \cite{HL71} implies that a hypersurface in the Riemannian
sphere $S^m$ is homogeneous if and only if it is a principal orbit
of the $s$-representation of an $(m+1)$-dimensional semisimple
Riemannian symmetric space $G/K$ of rank two. The classification
of homogeneous hypersurfaces in $S^m$ can therefore be easily
deduced from Cartan's classification of Riemannian symmetric
spaces.

If $G/K$ is Hermitian symmetric, then $m$ is odd, say $m = 2n+1$,
and the $s$-representation induces an action on the corresponding
complex projective space $\C P^n$ via the Hopf map $S^{2n+1} \to \C
P^n$. Takagi \cite{Ta73} showed in 1973 that a real hypersurface in
$\C P^n$ is homogeneous if and only if it is a principal orbit of an
action that is induced in this way from the $s$-representation of a
semisimple Hermitian symmetric space of rank two. Thus the
classification of homogeneous hypersurfaces in $\C P^n$ can easily
be deduced from the classification of Hermitian symmetric spaces.

A remarkable consequence of Takagi's result is that each
homogeneous hypersurface in $\C P^n$ is a Hopf hypersurface. A
real hypersurface $M$ of an almost Hermitian manifold $\bar{M}$ is
a Hopf hypersurface if the one-dimensional foliation on $M$
induced by the rank one distribution $J(\nu M)$ is totally
geodesic, where $\nu M$ is the normal bundle of $M$ and $J$ is the
complex structure on $\bar{M}$. In 1989 the first author
classified in \cite{Be89} the homogeneous Hopf hypersurfaces in
the complex hyperbolic space $\C H^n$. Any such Hopf hypersurface
is either a horosphere in $\C H^n$, or a tube around a totally
geodesic $\R H^n$ or $\C H^k$ for some $k \in \{0,\ldots,n-1\}$.

For some time it was believed that, as is the case for $\C P^n$,
every homogeneous hypersurface in $\C H^n$ was a Hopf hypersurface.
It came as a kind of surprise when Lohnherr \cite{Lo98}
constructed in 1998 a counterexample: the ruled real hypersurface
$W^{2n-1}$ in $\C H^n$ which is determined by a horocycle in a
totally geodesic $\R H^2 \subset \C H^n$ is a non-Hopf homogeneous
real hypersurface. Recently the first author and Tamaru
\cite{BT05} obtained the classification of homogeneous
hypersurfaces in $\C H^n$.
A connected real hypersurface in $\C H^n$, $n \geq 2$, is
homogeneous if and only if it is holomorphically congruent to
\begin{itemize}
\item[(A)] a tube around a totally geodesic $\CH[k]$ for some
$k\in\{0,\dots,n-1\}$, or

\item[(B)] a tube around a totally
geodesic $\R H^n$, or

\item[(H)] a horosphere in $\C H^n$, or

\item[(S)] the ruled real hypersurface $W^{2n-1}$ or one of its
equidistant hypersurfaces, or

\item[(W)] a tube around the minimal ruled submanifold
$W^{2n-k}_\varphi$ for some $\varphi \in (0,\pi/2]$ and $k \in
\{2,\ldots,n-1\}$, where $k$ is even if $\varphi \neq \pi/2$.
\end{itemize}

\noindent The construction of the minimal ruled submanifolds
$W^{2n-k}_{\varphi}$ will be described later in this article. We
just mention here that the normal bundle of $W^{2n-k}_\varphi$ has
rank $k$ and constant K\"{a}hler angle $\varphi$.

The hypersurfaces of type (A), (B) and (H) are Hopf hypersurfaces
and their geometry is well understood. The first author gave in
\cite{Be98} a Lie theoretic construction of the homogeneous
hypersurfaces of type (S) and investigated their geometry. The aim
of this paper is to investigate the geometry of the other
homogeneous hypersurfaces and their focal sets
$W^{2n-k}_{\varphi}$.

The motivation for our investigations originates from the
question: Is every real hypersurface with constant principal
curvatures in a complex hyperbolic space $\CH$ an open part of a
homogeneous hypersurface? \'{E}lie Cartan \cite{Ca38} gave an
affirmative answer for the corresponding question in real
hyperbolic space. Some time ago the first author demonstrated in
\cite{Be89} that for $\CH$ the answer is yes within the class of
Hopf hypersurfaces. For arbitrary real hypersurfaces, we obtained
recently an affirmative answer in \cite{BD06} in case of $\CH[2]$,
and for $n \geq 3$ we obtained an affirmative answer in
\cite{BD05} provided that the number $g$ of distinct principal
curvatures satisfies $g \leq 3$. It is a well-established fact
that the classification problem of hypersurfaces with constant
principal curvatures is intimately related to the understanding of
the geometric structure of their focal sets. Therefore the next
step is to understand more thoroughly the geometry of the
homogeneous hypersurfaces in $\CH$ and their focal sets. The two
main results of this paper are as follows. (1) We determine
explicitly the principal curvatures and their multiplicities for
all homogeneous hypersurfaces in $\CH$. A consequence is that $g
\in \{2,3,4,5\}$ for any such hypersurface. (2) We give a
characterization of the non-totally geodesic focal sets of
homogeneous hypersurfaces in $\CH$ in terms of their second
fundamental form.

We now describe the contents of this paper. In Section 2 we
summarize some basic material about the complex hyperbolic space.
In Section 3 we characterize the minimal ruled submanifolds
$W_\varphi^{2n-k}$ in terms of the second fundamental form. More
precisely, we prove

\begin{thm}[Rigidity of the submanifold $W_\varphi^{2n-k}$]
Let $M$ be a
$(2n-k)$-dimensional connected submanifold in $\CH$, $n\geq 2$, with
normal bundle $\nu M\subset T\CH$ of constant K\"{a}hler angle
$\varphi\in(0,\pi/2]$. Assume that there exists a unit vector field
$Z$ tangent to the maximal complex distribution on $M$ such that the
second fundamental form $\II$ of $M$ is given by the trivial
symmetric bilinear extension of
\[ 2\II(Z,P\xi)=\sin^2(\varphi)\,\xi \]
for all $\xi\in \nu M$,
where $P\xi$ is the tangential component of $J\xi$.
Then $M$ is holomorphically congruent to
an open part of the ruled minimal submanifold $W_\varphi^{2n-k}$.
Conversely, the second fundamental form of $W_\varphi^{2n-k}$ is
of this form.
\end{thm}

This shows that the minimal ruled submanifold $W_\varphi^{2n-k}$
has three distinct constant principal curvatures
$\sin(\varphi)/2$, $-\sin(\varphi)/2$ and $0$ with multiplicities
$1$, $1$ and $2n-k-2$ with respect to each unit normal vector. In
Section 4 we determine the principal curvatures of the tubes
around $W_\varphi^{2n-k}$ together with their multiplicities. A
table containing the principal curvatures and their multiplicities
of all homogeneous hypersurfaces in $\C H^n$ is given at the end
of the paper.

The investigation of the geometry of the tubes around $\W{2n-k}$ for
$0 < \varphi < \pi/2$ can be reduced to the one of the tubes around
$\W{4}$ in $\C H^3$. In the final part of the paper we derive some
geometric information about these particular  tubes in terms of some
autoparallel distributions of their tangent bundles.

The second author has been supported by project BFM 2003-02949
(Spain) and an EMBARK postdoctoral fellowship from the Irish
Research Council for Science, Engineering and Technology.

\section{Preliminaries}

In this section we summarize some basic facts about the
complex hyperbolic space. For details we refer to \cite{BTV95}.

Let $\CH$, $n \geq 2$, be the $n$-dimensional complex hyperbolic
space equipped with the Fubini-Study metric $\langle \cdot , \cdot
\rangle$ of constant holomorphic sectional curvature $-1$. The
Riemannian curvature tensor $\bar{R}$ on $\C H^n$ is given by
\[
4\bar{R}(X,Y)Z =
\langle X,Z \rangle Y -\langle Y,Z \rangle X
+ \langle JX,Z \rangle JY - \langle JY,Z \rangle JX
+ 2\langle JX,Y \rangle JZ,
\]
where $J$ is the complex structure on $\C H^n$.

We denote by $\CH(\infty)$ the ideal boundary of $\CH$. Each
element $x$ of $\CH(\infty)$ is an equivalence class of asymptotic
geodesics in $\CH$. We equip $\CH\cup\CH(\infty)$ with the cone
topology. Then $\CH\cup\CH(\infty)$ is homeomorphic to a closed
ball in the Euclidean space $\R^{2n}$. For each $o\in\CH$ and each
$x\in\CH(\infty)$ there exists a unique geodesic $\gamma_{ox} : \R
\to \CH$ such that $\lVert\dot{\gamma}_{ox}\rVert= 1$,
$\gamma_{ox}(0) = o$ and
$\lim_{t\rightarrow\infty}\gamma_{ox}(t)=x$.

The connected component of the isometry group of $\CH$ is the
special unitary group $G = SU(1,n)$. We fix a point $o \in \CH$
and denote by $K$ the isotropy subgroup of $G$ at $o$. Then $K$ is
isomorphic to $S(U(1)\times U(n)) \subset SU(1,n)$ and $(G,K)$ is
a symmetric pair. Let ${\mathfrak g} = {\mathfrak k} + {\mathfrak
p}$ be the corresponding Cartan decomposition of the Lie algebra
${\mathfrak g}$ of $G$. As usual we identify $T_o\CH$ with
${\mathfrak p}$.

We now fix a point $x \in \CH(\infty)$ and denote by ${\mathfrak a}$
the one-dimensional linear subspace of ${\mathfrak p}$ spanned by
$\dot{\gamma}_{ox}(0) \in T_o\CH \cong {\mathfrak p}$. As the rank
of $\CH$ is one, ${\mathfrak a}$ is a maximal abelian subspace of
${\mathfrak p}$. Let $ {\mathfrak g} = {\mathfrak g}_{-2\alpha} +
{\mathfrak g}_{-\alpha} + {\mathfrak g}_0 + {\mathfrak g}_\alpha +
{\mathfrak g}_{2\alpha}$ be the root space decomposition of
${\mathfrak g}$ induced by ${\mathfrak a}$. Then ${\mathfrak n} =
{\mathfrak g}_\alpha + {\mathfrak g}_{2\alpha}$ is a $2$-step
nilpotent subalgebra of ${\mathfrak g}$ which is isomorphic to the
$(2n-1)$-dimensional Heisenberg algebra. The center of ${\mathfrak
n}$ is the one-dimensional subalgebra ${\mathfrak g}_{2\alpha}$.
Moreover, ${\mathfrak g} = {\mathfrak k} + {\mathfrak a} +
{\mathfrak n}$ is an Iwasawa decomposition of ${\mathfrak g}$. We
denote by $A$ and $N$ the connected closed subgroup of $G$ with Lie
algebra ${\mathfrak a}$ and ${\mathfrak n}$, respectively. The orbit
$A \cdot o$ of $A$ through $o$ is just the path $\gamma_{ox}(\R)$ of
the geodesic $\gamma_{ox}$ in $\CH$, and the orbits of $N$ are the
horospheres in $\CH$ centered at $x$. The solvable subgroup $AN
\subset KAN = G$ acts simply transitively on $\CH$. Thus we can
identify ${\mathfrak a} + {\mathfrak n}$ with $T_o\CH$. The
Riemannian metric on $\CH$ induces an inner product $\langle \cdot ,
\cdot \rangle$ on ${\mathfrak a} + {\mathfrak n}$, and we may
identify $\CH$ with the solvable Lie group $AN$ equipped with the
left-invariant Riemannian metric which is induced from $\langle
\cdot , \cdot \rangle$.

The complex structure $J$ on $\CH$ induces a complex structure on
the vector space ${\mathfrak a} + {\mathfrak n} \cong T_o\CH$ which
we will also denote by $J$. We define $B:=\dot{\gamma}_{ox}(0) \in
{\mathfrak a}$ and $Z := JB \in {\mathfrak g}_{2\alpha}$. Note that
${\mathfrak g}_\alpha$ is $J$-invariant. The Lie algebra structure
on ${\mathfrak a} + {\mathfrak n}$ is given by the trivial
skew-symmetric bilinear extension to $({\mathfrak a} + {\mathfrak
n}) \times ({\mathfrak a} + {\mathfrak n})$ of the relations
\begin{equation}\label{bracket}
[B,Z] = Z\ ,\ 2[B,U] = U\ ,\ [U,V] = \langle JU,V \rangle Z
\qquad (U,V \in {\mathfrak g}_\alpha).
\end{equation}
This shows that ${\mathfrak a} + {\mathfrak n}$ is a semidirect sum
of the two Lie algebras ${\mathfrak a}$ and ${\mathfrak n}$. Let
$\Exp_{\mathfrak n} : {\mathfrak n} \to N$ be the Lie exponential
map. The group structure on the semidirect product $AN$ is given by
\begin{equation}\label{expsolvprod}
\begin{array}{l}
\Bigl(a,\Exp_{\g{n}}(U+xZ)\Bigr)\cdot
\Bigl(b,\Exp_{\g{n}}(V+yZ)\Bigr)\\
\displaystyle\qquad=\Bigl(a+b,\Exp_{\g{n}} \Bigl(U+e^{a/2}V +
\Bigl(x+e^ay+\frac{1}{2}e^{a/2}\langle JU,V \rangle\Bigr) Z
\Bigr)\Bigr)
\end{array}
\end{equation}
for all $a,b,x,y \in \R$ and $U,V \in {\mathfrak g}_\alpha$. Here we
identify the one-dimensional Lie group $A$ in the canonical way with
$\R$ such that $1 \in \R$ corresponds to $\Exp_{\mathfrak a}(B)$
with the Lie exponential map $\Exp_{\mathfrak a} : {\mathfrak a} \to
A$. Finally, the Lie exponential map $\Exp_{{\mathfrak a} +
{\mathfrak n}} : {\mathfrak a} + {\mathfrak n} \to AN$ is given by
\begin{equation}\label{expsolv}
\Exp_{{\mathfrak a} + {\mathfrak n}}\bigl(aB+U+xZ\bigr)
=\Bigl(a,\Exp_{\g{n}}\Bigl(\rho(a/2)U+\rho(a)xZ\Bigr)\Bigr)
\end{equation}
for all $a,x \in \R$ and $U \in {\mathfrak g}_\alpha$, where
the analytic function $\rho:\R\to\R$ is defined by
\[
\rho(s)=\left\{
\begin{array}{cll}
\displaystyle\frac{e^s-1}{s}    &,& \mbox{if }s\neq 0,\\
\noalign{\medskip} 1                               &,& \mbox{if
}s=0.
\end{array}\right.
\]
The Lie exponential map $\Exp_{{\mathfrak a} + {\mathfrak n}}$ is
a diffeomeorphism. If $V \in {\mathfrak g}_\alpha$ is a unit vector, then
the geodesic $\gamma$ in $\CH$ with $\gamma(0) = o$ and $\dot{\gamma}(0) = V$
is given by
\begin{equation} \label{geodesic}
\gamma(t) = (\ln\sech^2(t/2),
\Exp_{\mathfrak n}(2\tanh(t/2)V)),
\end{equation}
and its tangent vector field $\dot{\gamma}$ is given by
\begin{equation} \label{geodesictangentvector}
\dot{\gamma}(t) = -\tanh(t/2)B +\sech(t/2)V.
\end{equation}

We denote by $\bar\nabla$ the Levi-Civita covariant derivative of
$\CH$. The standard method for calculating the Levi-Civita covariant
derivative of a Lie group equipped with a left-invariant Riemannian
metric yields
\begin{equation}\label{LeviCivita}
\begin{array}{rcl}
\enabla_{aB+U+xZ}\bigl(bB+V+yZ\bigr)
    & = &\displaystyle\Bigl(\frac{1}{2}\langle U,V\rangle+x\,y\Bigr)B
    -\frac{1}{2}\bigl(bU+yJU+xJV\bigr)
    \\
\noalign{\medskip} &&
\displaystyle+\Bigl(\frac{1}{2}\,\langle JU,V\rangle-b\,x\Bigr)Z,
\end{array}
\end{equation}
where $a,b,x,y\in\R$ and $U,V\in {\mathfrak g}_\alpha$ and all
elements in ${\mathfrak a}+{\mathfrak n}$ are considered as
left-invariant vector fields on $AN \cong \CH$.

\section{The ruled submanifolds $W^{2n-k}$ and $W^{2n-k}_{\varphi}$}

The submanifolds $W^{2n-k}$ and $W^{2n-k}_{\varphi}$ were first
constructed by the first author and Br\"{u}ck in \cite{BB01}.

Let $\g{v}$ be a linear subspace of ${\mathfrak g}_{\alpha}$. For
each $0 \neq v\in \g{v}$ the K\"{a}hler angle of $\g{v}$ with
respect to $v$ is the angle $\varphi(v)\in[0,\pi/2]$ between
$\g{v}$ and the real span of $Jv$. Thus $\varphi(v)\in[0,\pi/2]$
is determined by requiring that $\cos(\varphi(v))\lVert v\rVert$
is the length of the orthogonal projection of $Jv$ onto $\g{v}$.
We say that $\g{v}$ has \emph{constant K\"{a}hler angle} $\varphi$
if $\varphi(v)=\varphi$ for all nonzero vectors $v\in \g{v}$. The
subspaces of ${\mathfrak g}_{\alpha}$ with constant K\"{a}hler
angle $\varphi = 0$ are precisely the complex subspaces of
${\mathfrak g}_{\alpha}$, and the subspaces of ${\mathfrak
g}_{\alpha}$ with constant K\"{a}hler angle $\varphi = \pi/2$ are
precisely the real subspaces of ${\mathfrak g}_{\alpha}$.

Let $\g{w}$ be a linear subspace of ${\mathfrak g}_{\alpha}$ such
that the orthogonal complement $\g{w}^\perp = {\mathfrak
g}_{\alpha}\ominus \g{w}$ of $\g{w}$ in ${\mathfrak g}_{\alpha}$ has
constant K\"{a}hler angle $\varphi \in [0,\pi/2]$. Then $\g{s}=\g{a} +
\g{w} + {\mathfrak g}_{2\alpha}$ is a subalgebra of $\g{a} + \g{n}$.
Denote by $S$ the connected closed subgroup of $AN$ with Lie algebra
$\g{s}$ and by $N_K^0(S)$ the identity component of the normalizer
of $S$ in $K$. Then $N_K^0(S)S \subset KAN$ acts on $\CH$ with
cohomogeneity one. The orbit $W^{2n-k}_{\varphi} = N_K^0(S)S \cdot o
= S \cdot o$ of this action containing the point $o$ is a
$(2n-k)$-dimensional submanifold of $\C H^n$, where $k = \dim
\g{w}^\perp$.

If $\varphi = 0$, that is, $\g{w}^\perp$ is a complex subspace
of ${\mathfrak g}_{\alpha}$, then $W^{2n-k}_0$ is a totally geodesic
complex hyperbolic subspace $\C H^{n-k'}$, where $k=2k'$.

If $\varphi = \pi/2$, then $\g{w}^\perp$ is a $k$-dimensional real
subspace of ${\mathfrak g}_{\alpha}$. If $k=1$, then
$W^{2n-1}_{\pi/2}$ is the ruled real hypersurface $W^{2n-1}$
determined by a horocycle in a totally geodesic $\R H^2 \subset \C
H^n$. The orbits of $N_K^0(S)S$ form a homogeneous codimension one
foliation on $\C H^n$ whose geometry has been investigated by the
first author in \cite{Be98}. If $k > 1$, then $W^{2n-k}_{\pi/2}$ is
a $(2n-k)$-dimensional homogeneous submanifold of $\C H^n$ with
totally real normal bundle of rank $k$. We will sometimes use the
notation $W^{2n-k} := W^{2n-k}_{\pi/2}$.

If $0 < \varphi < \pi/2$, then $k$ is even and $W^{2n-k}_{\varphi}$
is a $(2n-k)$-dimensional homogeneous submanifold of $\C H^n$ whose
normal bundle has constant K\"{a}hler angle $\varphi$ and rank $k$.

As $\C H^n$ is a two-point homogeneous space, the construction of
the submanifolds $W^{2n-k}$ and $W^{2n-k}_{\varphi}$ does not depend
on the choice of the two points $o \in \C H^n$ and $x \in \C
H^n(\infty)$, or equivalently, on the choice of the Iwasawa
decomposition of $G$. Our next aim is to investigate the geometry of
the submanifolds $W^{2n-k}_{\varphi}$, $0 < \varphi \leq \pi/2$.

Let $\C\g{w}^\perp$ be the complex subspace of ${\mathfrak
g}_{\alpha}$ spanned by $\g{w}^\perp$ and ${\mathfrak d} =
\C\g{w}^\perp \ominus \g{w}^\perp$ be the orthogonal complement of
$\g{w}^\perp$ in $\C\g{w}^\perp$. As $\varphi > 0$, we have $k =
\dim_{\C}\C\g{w}^\perp$ and hence $k = \dim {\mathfrak w}^\perp =
\dim {\mathfrak d}$. For each $\xi \in {\mathfrak w}^\perp$ we
decompose $J\xi \in \C{\mathfrak w}^\perp = {\mathfrak d} +
{\mathfrak w}^\perp$ into $J\xi = P\xi + F\xi$ with $P\xi \in
{\mathfrak d}$ and $F\xi \in {\mathfrak w}^\perp$. Since ${\mathfrak
w}^\perp$ has constant K\"{a}hler angle $\varphi$, we have $\langle F\xi
, F\xi \rangle = \cos^2(\varphi) \langle \xi , \xi \rangle$ and
hence $\langle P\xi , P\xi \rangle = \sin^2(\varphi) \langle \xi ,
\xi \rangle$. As $\varphi > 0$, the homomorphism $P : {\mathfrak
w}^\perp \to {\mathfrak d}$ is injective, and as $\dim {\mathfrak
w}^\perp = \dim {\mathfrak d}$ we see that $P : {\mathfrak w}^\perp
\to {\mathfrak d}$ is an isomorphism. From $-\xi = JJ\xi = JP\xi +
JF\xi = JP\xi + PF\xi + F^2\xi$ we see that the ${\mathfrak
d}$-component $(JP\xi)_{\mathfrak d}$ of $JP\xi$ is equal to
$-PF\xi$, and hence $\langle (JP\xi)_{\mathfrak d} ,
(JP\xi)_{\mathfrak d} \rangle = \langle PF\xi , PF\xi \rangle =
\sin^2(\varphi)\langle F\xi,F\xi \rangle =
\sin^2(\varphi)\cos^2(\varphi)\langle \xi,\xi \rangle =
\cos^2(\varphi) \langle P\xi , P\xi \rangle$. Since $P : {\mathfrak
w}^\perp \to {\mathfrak d}$ is an isomorphism, this implies that
${\mathfrak d}$ has constant K\"{a}hler angle $\varphi$ as well.

We denote by ${\mathfrak c}$ the maximal complex subspace of
${\mathfrak s}$. Note that ${\mathfrak a} + {\mathfrak g}_{2\alpha}
\subset {\mathfrak c}$, $\dim_{\C}{\mathfrak c} = n-k$ and
${\mathfrak s} = {\mathfrak c} + {\mathfrak d}$. Then we have the
orthogonal decomposition
\[ {\mathfrak a} + {\mathfrak n} = {\mathfrak c} + {\mathfrak d}
+ {\mathfrak w}^\perp.\] We denote by ${\mathfrak A}$, ${\mathfrak
C}$, ${\mathfrak D}$ and ${\mathfrak W}^\perp$ the left-invariant
distributions on $\CH$ along $\W{2n-k}$ which are induced by
${\mathfrak a}$, ${\mathfrak c}$, ${\mathfrak d}$ and ${\mathfrak
w}^\perp$, respectively. By construction, we have ${\mathfrak C} +
{\mathfrak D} = T\W{2n-k}$ and ${\mathfrak W}^\perp=\nu \W{2n-k}$.

\begin{proposition}\label{thCharruled}
The submanifold $\W{2n-k}$, $0 < \varphi \leq \pi/2$, of $\CH$ has
the following properties:
\begin{itemize}
\item[(i)] The maximal holomorphic subbundle $\g{C}$ of
$T\W{2n-k}$ is autoparallel and the leaves of the induced foliation
on $\W{2n-k}$ are totally geodesic $\CH[n-k] \subset \CH$. Hence
$\W{2n-k}$ is a ruled submanifold of $\CH$.

\item[(ii)] The following statements are equivalent:
\begin{itemize}
\item[(a)] the distribution $\g{D}$ on $\W{2n-k}$ is
integrable;

\item[(b)] the distribution $\g{A} + \g{D}$ on
$\W{2n-k}$ is integrable;

\item[(c)] the normal bundle $\g{W}^\perp$ is flat with respect to
the normal connection;

\item[(d)] $\varphi = \pi/2$.
\end{itemize}
In this case the leaves of the foliation on $\W[\pi/2]{2n-k}$
induced by $\g{A} + \g{D}$ are totally geodesic $\R H^{k+1} \subset
\CH$ and the leaves of the foliation on $\W[\pi/2]{2n-k}$ induced by
$\g{D}$ are horospheres with center $x$ in these totally geodesic
$\R H^{k+1} \subset \CH$.

\item[(iii)] For each $0 \neq \xi \in \g{w}^\perp$ the
left-invariant distribution ${\mathfrak A} + \R P\xi$ on $\W{2n-k}$
is autoparallel and the leaves of the induced foliation on
$\W{2n-k}$ are totally geodesic $\R H^2 \subset \CH$.

\item[(iv)] For each $0 \neq \xi \in \g{w}^\perp$ the
left-invariant distribution $\R P\xi$ on $\W{2n-k}$ is integrable
and the leaves of the induced foliation on $\W{2n-k}$ are
horocycles with center $x$ in the totally geodesic $\R H^2 \subset
\CH$ given by the distribution ${\mathfrak A} + \R P\xi$.
\end{itemize}
\end{proposition}

\begin{proof}
(i) follows immediately from (\ref{LeviCivita}) and the fact that
the only complex totally geodesic submanifolds of $\CH$ are
complex hyperbolic spaces.

From (\ref{bracket}) we see that $2[aB+U,bB+V] = aV - bU + 2\langle
JU,V \rangle Z$ for all $a,b \in \R$ and $U,V \in {\mathfrak D}$.
This shows that ${\mathfrak A} + \g{D}$ is integrable if and only if
$\g{D}$ is integrable if and only if $\g{D}$ is real, that is,
$\varphi = \pi/2$. On the other hand, \eqref{LeviCivita} implies
$2\nabla^\perp_{aB+U+xZ}\xi=-x F\xi$ for all $a,x\in\R$, $U\in\g{D}$
and $\xi\in\g{W}^\perp$. Hence, the normal bundle $\g{W}^\perp$ is
flat if and only if $F=0$, or equivalently, $\varphi=\pi/2$. In this
case (\ref{LeviCivita}) yields
\[
2\enabla_{aB+U}(bB+V) = \langle U,V \rangle B - bU \in {\mathfrak
A}+ \g{D}
\]
for all $a,b \in \R$ and $U,V \in {\mathfrak D}$. This shows that
${\mathfrak A}+ \g{D}$ is autoparallel and its leaves are totally
geodesic real submanifolds of $\CH$. The only real totally geodesic
submanifolds of $\CH$ are real hyperbolic spaces. Finally, for all
$U,V \in \g{D}$ we have $2\enabla_UV = \langle U,V \rangle B$ and
$2\enabla_UB = -U$, which implies that the leaves of $\g{D}$ are
spherical hypersurfaces of the corresponding real hyperbolic
subspaces. Since the sectional curvature of a totally geodesic real
hyperbolic subspace is $-1/4$, and the mean curvature vector field
of any leaf of $\g{D}$ is $(1/2)B$, it follows that the leaves of
$\g{D}$ are horospheres centered at $x$ in the real hyperbolic
subspaces. This finishes the proof of (ii).

For any $aB+xP\xi,bB+yP\xi \in {\mathfrak A} + \R P\xi$ we get from
(\ref{LeviCivita}) and using
$\langle P\xi , P\xi \rangle
= \sin^2(\varphi) \langle \xi , \xi \rangle$ that
$2\enabla_{aB+xP\xi}(bB+yP\xi)
= xy\sin^2(\varphi)B - bxP\xi \in {\mathfrak A} + \R P \xi$.
From this we easily get the assertion (iii).

Finally, define $U_\xi = P\xi/\sin(\varphi)$. Then
(\ref{LeviCivita}) implies $2\enabla_{U_\xi}U_\xi = B$ and
$4\enabla_{U_\xi}\enabla_{U_\xi}U_\xi = -U_\xi$. Since the real
hyperbolic planes in (iii) have constant sectional curvature
$-1/4$, this shows that the integral curves of $U_\xi$ are
horocycles with center $x$ in the corresponding real hyperbolic
planes. This proves (iv).
\end{proof}

The previous proposition implies a nice geometric construction of
the ruled submanifolds $\W{2n-k} \subset \CH$.

\begin{corollary}\label{thConstruction}
Let $k \in \{1,\ldots,n-1\}$, and fix a totally geodesic $\CH[n-k]
\subset \CH$ and points $o \in \CH[n-k]$ and $x \in
\CH[n-k](\infty)$. Let $KAN$ be the Iwasawa decomposition of
$SU(1,n)$ with respect to $o$ and $x$, and let $H'$ be the
subgroup of $AN$ which acts simply transitively on $\CH[n-k]$.
Next, let $W$ be a subspace of $\nu_o\CH[n-k]$ with constant
K\"{a}hler angle $\varphi \in (0,\pi/2]$ such that $\C W =
\nu_o\CH[n-k]$. Left translation of $W$ by $H'$ to all points in
$\CH[n-k]$ determines a subbundle ${\mathfrak V}$ of the normal
bundle $\nu \CH[n-k]$. At each point $p \in \CH[n-k]$ attach the
horocycles determined by $x$ and the linear lines in ${\mathfrak
V}_p$. The resulting subset $M$ of $\CH$ is holomorphically
congruent to the ruled submanifold $\W{2n-k}$.
\end{corollary}

\begin{proof}
Let $\W{2n-k}$ be the ruled minimal submanifold of $\CH$ constructed
from the Iwasawa decomposition $KAN$ associated with $o$ and $x$ and
from the choice of $\g{w}^\perp = \nu_o\CH[n-k] \ominus W$. We use
the above notations. We will show that $M = \W{2n-k}$. Let $p \in
\W{2n-k}$. Then there exists an isometry $s \in S$ with $p = s(o)$.
There is a unique vector $X$ in the Lie algebra $\g{s}$ of $S$ such
that $s = \Exp_{\g{a}+\g{n}}(X)$. We can write $X = aB + U + V + zZ$
with some $U \in \g{c} \ominus ({\mathfrak a} + {\mathfrak
g}_{2\alpha})$, $V \in \g{d}$, and $a,z \in \R$. Note that $[V,U]
=0$ because they are complex orthogonal. We now define
$g=\Exp_{\g{a}+\g{n}}(\rho(a/2)V)$
and $h=\Exp_{\g{a}+\g{n}}(aB+U+zZ)$.
Note that $h \in H'$. Using (\ref{expsolvprod}) and
(\ref{expsolv}) we get
\begin{eqnarray*}
gh &=& \Exp_{\g{a}+\g{n}}\left(\rho\left(\frac{a}{2}\right)V\right)
    \Exp_{\g{a}+\g{n}}(aB+U+zZ)\\
\noalign{\smallskip} &=& \left(0,\Exp_{\g{n}}
    \left(\rho\left(\frac{a}{2}\right)V\right)\right)
    \cdot\left(a,\Exp_{\g{n}}\left(\rho\left(\frac{a}{2}\right)U
    +\rho(a)z\,Z
    \right)\right)\\
\noalign{\smallskip} &=& \left(a,\Exp_{\g{n}}\left(
    \rho\left(\frac{a}{2}\right)(U+V)+ \rho(a)z\,Z\right)\right)\\
\noalign{\medskip} &=& \Exp_{\g{a}+\g{n}}(aA+U+V+zZ)\
    =\ \Exp_{\g{a}+\g{n}}(X)\ =\ s.
\end{eqnarray*}
By construction, $h(o) \in \CH[n-k]$, and $s(o) = g(h(o))$ is on the
horocycle with center $x$ through $h(o)$ tangent to $\R V$. From this
we conclude that $\W{2n-k} \subset M$. From Proposition
\ref{thCharruled} we already know that $M \subset \W{2n-k}$.
Altogether this implies $M = \W{2n-k}$.
\end{proof}

We now describe the geometry of $\W{2n-k}$ in terms of the second
fundamental form.

\begin{proposition} \label{secondfunform}
The second fundamental form $\II$ of $\W{2n-k}$ is given by
\[
\II\bigl(aB+U+P\xi+xZ,bB+V+P\eta+yZ\bigr)
=\frac{\sin^2(\varphi)}{2}\bigl(y\xi+x\eta\bigr)
\]
for all $\xi,\eta\in \g{w}^\perp$, $U,V \in \g{c}\ominus ({\mathfrak
a} + {\mathfrak g}_{2\alpha})$ and $a,b,x,y\in\R$. Thus $\II$ is
given by the trivial symmetric bilinear extension of
$2\II(Z,P\xi)=\sin^2(\varphi)\xi$ for all $\xi\in \g{w}^\perp$.
\end{proposition}

\begin{proof}
We denote by $(\cdot)^\perp$ the orthogonal projection onto
$\nu\W{2n-k}$. The Gau\ss\ formula and (\ref{LeviCivita}) imply
\begin{eqnarray*}
\II\bigl(aB+U+P\xi+xZ,bB+V+P\eta+yZ\bigr)
&=& \left(\enabla_{aB+U+P\xi+xZ}(bB+V+P\eta+yZ)\right)^\perp\\
&=& -\left(\frac{y}{2}JP\xi+\frac{x}{2}JP\eta\right)^\perp\\
&=& \frac{\sin^2(\varphi)}{2}\bigl(y\xi+x\eta\bigr),
\end{eqnarray*}
since $(JP\xi)^\perp =-\sin^2(\varphi)\xi$, which follows from the
fact that $\g{w}^\perp$ has constant K\"{a}hler angle $\varphi$.
\end{proof}

As an immediate consequence we get

\begin{corollary}
$W_\varphi^{2n-k}$ is a minimal ruled submanifold of $\CH$.
\end{corollary}

For $k > 1$ the previous corollary follows also from the general
fact that each singular orbit of a cohomogeneity one action is a
minimal submanifold. We will now show that the above equation for
the second fundamental form in fact characterizes the minimal ruled
submanifolds $W_\varphi^{2n-k}$ in $\CH$.

\begin{theorem}[Rigidity of the submanifold $W_\varphi^{2n-k}$]
\label{thRigidity}
Let $M$ be a $(2n-k)$-dimensional connected submanifold in $\CH$,
$n\geq 2$, with normal bundle $\nu M\subset T\CH$ of constant K\"{a}hler
angle $\varphi\in(0,\pi/2]$. Assume that there exists a unit vector
field $Z$ tangent to the maximal complex distribution on $M$ such
that the second fundamental form $\II$ of $M$ is given by the
trivial symmetric bilinear extension of
\begin{equation}\label{II}
2\II(Z,P\xi)=\sin^2(\varphi)\,\xi
\end{equation}
for all $\xi\in \nu M$,
where $P\xi$ is the tangential component of $J\xi$.
Then $M$ is holomorphically congruent to
an open part of the ruled minimal submanifold $W_\varphi^{2n-k}$.
Conversely, the second fundamental form of $W_\varphi^{2n-k}$ is
of this form.
\end{theorem}

\begin{proof}
The last statement is a consequence of Proposition
\ref{secondfunform}. For the other part we use Corollary
\ref{thConstruction}.

We decompose the tangent bundle $TM$ of $M$ orthogonally into $TM =
{\mathfrak C} + {\mathfrak D}$, where ${\mathfrak C}$ is the maximal
complex subbundle of $TM$. For each $\xi \in \Gamma(\nu M)$ we
decompose $J\xi$ orthogonally into $J\xi = P\xi + F\xi$ with $P\xi
\in \Gamma({\mathfrak D})$ and $F\xi \in \Gamma(\nu M)$. As above
one can show that ${\mathfrak D}$ has constant K\"{a}hler angle
$\varphi$ as well, and the bundle homomorphisms $P : \nu M \to
{\mathfrak D}$ and $F : \nu M \to \nu M$ are homomorphisms
satisfying $\langle P\xi , P\xi \rangle = \sin^2(\varphi)\langle \xi
, \xi \rangle$ and $\langle F\xi , F\xi \rangle=
\cos^2(\varphi)\langle \xi , \xi \rangle$. Note that if
$\varphi=\pi/2$ then $F$ is trivial and $P=J|_{\nu M}$; otherwise,
$P$ and $F$ are isomorphisms.

For all $U,V \in \Gamma(\g{C})$ and $\xi \in \Gamma(\nu M)$ we have,
using \eqref{II} and $\enabla J=0$,
\[
\langle\enabla_UV,\xi\rangle=\langle\II(U,V),\xi\rangle = 0
\]
and
\[
\langle\enabla_UV,J\xi\rangle = -\langle J\bar{\nabla}_UV,\xi\rangle
= - \langle \II(U,JV),\xi \rangle = 0.
\]
This shows that $\g{C}$ is an autoparallel subbundle of $TM$ and
each integral manifold is a totally geodesic submanifold of $\CH$.
As $\g{C}$ is a complex subbundle of complex rank $n-k$, each of
these integral manifolds must be an open part of a totally geodesic
$\CH[n-k]\subset \CH$.

Let $o \in M$ and $\mathcal{F}_o$ be the leaf of $\g{C}$ through
$o$, which is an open part of a totally geodesic $\CH[n-k] \subset
\CH$. Let $\gamma : I \to\mathcal{F}_o$ be a curve with $\gamma(0) =
o$. We prove that the normal spaces of $M$ along $\gamma$ are
uniquely determined by the differential equation
\begin{equation}\label{diffeq}
2\enabla_{\dot\gamma}X +\langle \dot\gamma , Z \rangle JX = 0
\end{equation}
along $\gamma^*\nu\mathcal{F}_o$. Let $X\in\Gamma(TM)$ and
$\xi\in\Gamma(\nu M)$. Using (\ref{II}) we get
\[
-\langle
\bar{\nabla}_{\dot{\gamma}}\xi,X \rangle = \langle
\II(\dot{\gamma},X),\xi \rangle = \langle \dot{\gamma},Z \rangle
\langle X,P\xi \rangle \langle \II(Z,P\xi),\xi \rangle/\langle
P\xi,P\xi\rangle = \frac{1}{2}\langle \dot{\gamma},Z \rangle \langle
P\xi,X \rangle,
\]
which implies
\begin{equation}\label{eqDiffEq2}
\enabla_{\dot{\gamma}}\xi = -\frac{1}{2}\langle \dot{\gamma},Z
\rangle P\xi + \nabla_{\dot{\gamma}}^\perp \xi,
\end{equation}
where $\nabla^\perp$ is the normal connection of $M$. Now, let $X$
be a vector field along $\gamma$ with $X_0\in\nu_o M$ and satisfying
\eqref{diffeq}. We may write $X=U+J\eta+\xi$ with
$U\in\Gamma(\gamma^*\g{C})$, $\xi,\eta\in\Gamma(\gamma^*\nu M)$ and
$U_0=\eta_0=0$. Then, using \eqref{eqDiffEq2} and $\enabla J=0$, we
get
\begin{eqnarray*}
0
&=& 2\enabla_{\dot\gamma}X+\langle\dot\gamma,Z\rangle JX\\
&=& 2\enabla_{\dot\gamma}U
+2J\enabla_{\dot\gamma}\eta+2\enabla_{\dot\gamma}\xi
+\langle\dot\gamma,Z\rangle JU
+\langle\dot\gamma,Z\rangle J^2\eta+\langle\dot\gamma,Z\rangle J\xi\\
&=& 2\enabla_{\dot\gamma}U+\langle\dot\gamma,Z\rangle JU
+P\left(2\nabla_{\dot\gamma}^\perp\eta
    +\langle\dot\gamma,Z\rangle F\eta\right)\\
&&  +2\nabla_{\dot\gamma}^\perp\xi+\langle\dot\gamma,Z\rangle F\xi
+F\left(2\nabla_{\dot\gamma}^\perp\eta
    +\langle\dot\gamma,Z\rangle F\eta\right).
\end{eqnarray*}
We have that $2\enabla_{\dot\gamma}U+\langle\dot\gamma,Z\rangle JU$
is tangent to $\g{C}$ because $\g{C}$ is a complex autoparallel
distribution. Hence, it follows that
$2\enabla_{\dot\gamma}U+\langle\dot\gamma,Z\rangle JU=0$. Since
$U_0=0$, the uniqueness of solutions to ordinary differential
equations implies $U_{t}=0$ for all $t$ and thus $X$ is
normal to $\mathcal{F}_o$ along $\gamma$. Similarly, the component
tangent to $P\nu M$ yields $2\nabla_{\dot{\gamma}}^\perp \eta +
\langle \dot{\gamma},Z \rangle F\eta = 0$ and since
$\eta_{0} = 0$ we have $\eta_{t}=0$ for all $t$.
Hence, $X_{t}\in \nu_{\gamma(t)}M$ for any $t$, which proves
our previous assertion.

We define $B:=-JZ$. The tangent vector $B_o$ determines a point $x
\in \C H^n(\infty)$, and thus, $o$ and $x$ give rise to an Iwasawa
decomposition of the Lie algebra of the isometry group of $\CH$.
Let us consider the ruled submanifold $\W{2n-k}$ determined by the
above Iwasawa decomposition and the normal space
$\g{w}^\perp=\nu_o M$ of constant K\"ahler angle $\varphi$. The
leaf $\mathcal{F}_o$ is an open part of the totally geodesic
$\CH[n-k]$ tangent to the maximal complex distribution of
$\W{2n-k}$ at $o$. We have just proved that \eqref{diffeq}
determines the normal bundle of a submanifold satisfying all the
hypotheses of Theorem \ref{thRigidity}, which implies $\nu_p
M=\nu_p\W{2n-k}$ for all $p\in\mathcal{F}_o$, that is, $\nu_p M$
is obtained by left translation of $\nu_o M$ for all
$p\in\mathcal{F}_o$. According to Corollary \ref{thConstruction}
it just remains to prove that at any point $p\in\mathcal{F}_o$ the
horocycles determined by $x$ and the linear lines in $P\nu_pM$ are
contained in a neighborhood of $p$ in $M$.

We now prove that the vector field $B$ is a geodesic vector field
and all its integral curves are geodesics in $\CH$ converging to
the point $x\in \C H^n(\infty)$. Since $B$ belongs to the maximal
complex distribution we have $\enabla_BB\in\Gamma(\g{C})$. Since
$B$ is a unit vector $\langle\enabla_BB,B\rangle=0$. Let $X \in
\Gamma({\mathfrak C} \ominus {\mathbb R}B)$ and $\eta \in
\Gamma(\nu M)$ be a local unit normal vector field of $M$. Using
the explicit expression for $\bar{R}$, the Codazzi equation,
(\ref{II}) and $\enabla J = 0$ we get
\begin{eqnarray*}
0 &=& 2\eR_{B P\eta JX \eta} \ =\ 2\langle
(\nabla_B^\perp\II)(P\eta,JX)-(\nabla_{P\eta}^\perp\II)(B,JX),\eta
\rangle \\
& = & -2\langle \II (P\eta,\nabla_BJX),\eta \rangle
\ = \ -2\langle\nabla_BJX,Z\rangle\langle\II(P\eta,Z),\eta\rangle\\
& = & -\sin^2(\varphi)\langle \enabla_BJX , Z \rangle \ =\
\sin^2(\varphi)\langle \enabla_BB , X \rangle.
\end{eqnarray*}
This implies $\langle\enabla_BB,X\rangle=0$ and therefore
$\enabla_BB=0$, which means that the integral curves of $B$ are
geodesics in $\CH$.

Now let $X\in\Gamma(TM\ominus\R B)$ and $\gamma$ an integral curve
of $X$. We consider the geodesic variation
$F(s,t)=\exp_{\gamma(s)}(tB_{\gamma(s)})$ of $\alpha(t)=F(0,t)$, the
geodesic in $\CH$ with initial condition $\dot\alpha(0)=B_o$. We
prove that $d(\alpha(t),F(s,t))$ tends to $0$ as $t$ goes to
infinity, where $d$ is the Riemannian distance function of $\CH$.

The transversal vector field of the geodesic variation $F$,
$\zeta(t)=(\partial F/\partial s)(0,t)$, is a Jacobi field along
$\alpha$ (thus, $4\zeta''-\zeta-3\langle \zeta,Z\rangle Z=0$)
with initial conditions $\zeta(0)=X_{\gamma(0)}$,
$\zeta'(0)=\enabla_{X_{\gamma(0)}} B$. Hence, we need to calculate
$\enabla_XB$.

Let $\eta\in\Gamma(\nu M)$ be a local unit vector field. Using
\eqref{II} we have
$\langle\enabla_XB,\eta\rangle=\langle\II(X,B),\eta\rangle=0$. We
also have $\langle\enabla_XB,B\rangle=0$. Now, using \eqref{II} we
get
\begin{equation}\label{eqDXBPxi}
\everymath{\displaystyle}
\begin{array}{rcl}
2\langle\enabla_XB,P\eta\rangle
&=& -2\langle\enabla_{X}JZ,J\eta-F\eta\rangle
    \ =\ -2\langle\II(X,Z),\eta\rangle
            -2\langle\II(X,B),F\eta\rangle\\
&=& -2\langle X,P\eta\rangle
        \langle\II(P\eta,Z),\eta\rangle/\sin^2(\varphi)
    =-\langle X,P\eta\rangle.
\end{array}
\end{equation}
Next, let $Y\in\Gamma(\g{C}\ominus\R B)$ and assume that
$X\in\Gamma(\g{C}\ominus\R B)$. Given $\xi\in\Gamma(\nu M)$ we have
$\langle\nabla_{P\eta}JY,P\xi\rangle=\langle\II(P\eta,Y),\xi\rangle
-\langle\II(P\eta,JY),F\xi\rangle=\frac{1}{2}\langle
Y,Z\rangle\langle P\eta,P\xi\rangle$. This, the expression for
$\bar{R}$, the Codazzi equation, (\ref{II}) and $\bar\nabla J = 0$
imply
\begin{eqnarray*}
-\sin^2(\varphi)\langle X,Y\rangle
&=& 4\eR_{X P\eta JY \eta}
    \ =\ 4\langle(\nabla_X^\perp\II)(P\eta,JY)
        -(\nabla_{P\eta}^\perp\II)(X,JY),\eta\rangle\\
&=& -4\langle\II(P\eta,\nabla_XJY),\eta\rangle
    +4\langle\II(X,\nabla_{P\eta}JY),\eta\rangle\\
&=& -4\langle\nabla_XJY,Z\rangle\langle\II(P\eta,Z),\eta\rangle
    +4\langle X,Z\rangle\langle\II(Z,\nabla_{P\eta}JY),\eta\rangle\\
&=& 2\sin^2(\varphi)\langle\enabla_XB,Y\rangle
    +\sin^2(\varphi)\langle X,Z\rangle\langle Z,Y\rangle.
\end{eqnarray*}
Hence, if $X\in\Gamma(\g{C}\ominus\R B)$ we have, as
$\enabla_XB\in\Gamma(\g{C})$, that $2\enabla_XB=-X-\langle
X,Z\rangle Z$.

Now assume that $X\in\Gamma(P\nu M)$ and write $X=P\xi$ with
$\xi\in\Gamma(\nu M)$. We have $\langle \nabla_{JY}P\xi , Z \rangle
= - \langle \bar\nabla_{JY}Z , J\xi - F\xi \rangle = - \langle
\II(JY,B),\xi \rangle + \langle \II(JY,Z),F\xi \rangle=0$. This,
together with the explicit expression for $\bar{R}$, the Codazzi
equation, (\ref{II}) and $\bar\nabla J = 0$ implies
\begin{eqnarray*}
0 & = & 2\eR_{P\xi JY P\xi \xi}\ =\ 2\langle
(\nabla_{P\xi}^\perp\II)(JY,P\xi)
-(\nabla_{JY}^\perp\II)(P\xi,P\xi),\xi \rangle \\
& = & -2\langle \II (\nabla_{P\xi} JY,P\xi),\xi \rangle
+ 4 \langle \II(\nabla_{JY}P\xi,P\xi),\xi\rangle \\
& = & -2\langle \nabla_{P\xi} JY , Z \rangle \langle \II
(Z,P\xi),\xi \rangle + 4 \langle \nabla_{JY}P\xi , Z \rangle
\langle \II(Z,P\xi),\xi\rangle \\
& = & -\sin^2(\varphi)\langle \enabla_{P\xi}JY , Z \rangle \ = \
\sin^2(\varphi)\langle \enabla_{P\xi}B , Y \rangle.
\end{eqnarray*}
Thus we get $\langle \enabla_{P\xi}B,Y\rangle=0$, and as a
consequence, we get $2\enabla_{P\xi}B=-P\xi$ using \eqref{eqDXBPxi}.

All in all, this implies
\begin{equation}\label{eqDXB}
\enabla_XB=-\frac{1}{2}X-\frac{1}{2}\langle X,Z\rangle Z \quad
\mbox{for all $X\in\Gamma(TM\ominus\R B)$.}
\end{equation}

Therefore, if $X\in T_{\alpha(0)}M\ominus\R B_{\alpha(0)}$ is a unit
vector and $\mathcal{B}_X$ denotes $\bar\nabla$-parallel translation
of $X$ along $\alpha$, we get
\[
\zeta(t)=e^{-t/2}\mathcal{B}_X(t)
+(e^{-t}-e^{-t/2})\langle X,Z_{\alpha(0)}\rangle Z_{\alpha(t)}.
\]
Note that $Z_{\alpha(t)}$ is a parallel vector field along $\alpha$
since $Z_{\alpha(t)}'=
\enabla_{B_{\alpha(t)}}Z=J\enabla_{B_{\alpha(t)}}B=0$. We easily see
that $\lim_{t\to\infty}\lVert \zeta(t) \rVert = 0$, which implies
$\lim_{t \to \infty} d(\alpha(t),F(s,t)) = 0$ as
$d(\alpha(t),F(s,t)) \leq s\lVert \zeta(t) \rVert$. Altogether this
shows that the integral curves of $B$ are asymptotic geodesics
corresponding to the point $x\in \C H^n(\infty)$.

Now let $p\in\mathcal{F}_o$ and $\xi_p\in\nu_pM$ be a unit vector.
The theorem now follows if we prove that the horocycle determined
by $P\xi_p/\sin(\varphi)$ and the point $x\in \C H^n(\infty)$ is
locally contained in $M$. To achieve this we will construct a unit
local vector field $\xi\in\Gamma(\nu M)$ such that the previous
horocycle is an integral curve of $P\xi/\sin(\varphi)$.

Let $\gamma:I\to M$ be a curve in $M$ satisfying the differential
equation
\begin{equation}\label{eqGamma}
\nabla_{\dot\gamma}\dot\gamma
=\frac{1}{2}\langle\dot\gamma,\dot\gamma\rangle B,\qquad
\dot\gamma(0)=P\xi_p/\sin(\varphi).
\end{equation}
We first prove that $\gamma$ is parametrized by arc length and that
it remains tangent to $P\nu M$.

Write $\dot\gamma=aB+xZ+X+P\eta$ for some differentiable functions
$a,x:I\to\R$ , and vector fields
$X\in\Gamma(\gamma^*(\g{C}\ominus(\R B+\R Z)))$ and
$\eta\in\Gamma(\gamma^*\nu M)$. Since $Z=JB$, the definition of
$\gamma$ and \eqref{eqDXB} show
\[
\frac{dx}{dt}
=\frac{d}{dt}\langle\dot\gamma,Z\rangle
=\langle\nabla_{\dot\gamma}\dot\gamma,Z\rangle
    +\langle\nabla_{\dot\gamma}Z,\dot\gamma\rangle
=\langle xB-\frac{1}{2}JX-\frac{1}{2}JP\eta,\dot\gamma\rangle
=ax.
\]
Since $x(0)=0$, the uniqueness of solutions to ordinary differential
equations implies that $x(t)=0$ for all $t$.

Let $Y\in\Gamma(\R B+P\nu M)$ and $\zeta\in\Gamma(\nu M)$. Then,
\eqref{II} yields
$\langle\enabla_YX,\zeta\rangle=\langle\II(Y,X),\zeta\rangle=0$ and
$\langle\enabla_YX,J\zeta\rangle=-\langle\II(Y,JX),\zeta\rangle=0$.
Moreover, since $\enabla_YB\in\Gamma(P\nu M)$ by \eqref{eqDXB} we have
$\langle\enabla_YX,B\rangle=-\langle\enabla_YB,X\rangle=0$. Also,
$2\langle\enabla_XX,B\rangle=-2\langle\enabla_XB,X\rangle =\langle
X,X\rangle$ and
$\langle\enabla_XX,P\eta\rangle=-\langle\II(X,JX),\eta\rangle
-\langle\II(X,X),F\eta\rangle=0$. Hence,
\begin{eqnarray*}
\frac{d}{dt}\langle X,X\rangle
&=& \frac{d}{dt}\langle\dot\gamma,X\rangle
\ =\ \langle\nabla_{\dot\gamma}\dot\gamma,X\rangle
    +\langle\nabla_{\dot\gamma}X,\dot\gamma\rangle\\
&=& a\langle\enabla_{\dot\gamma}X,B\rangle
    +\langle\enabla_{\dot\gamma}X,X\rangle
    +\langle\enabla_{\dot\gamma}X,P\eta\rangle
\ =\ \langle\enabla_{\dot\gamma}X,X\rangle
    +a\langle\enabla_{X}X,B\rangle\\
&=& \langle\nabla_{\dot\gamma}X,X\rangle
    +\frac{a}{2}\langle X,X\rangle
\ =\ \frac{1}{2}\frac{d}{dt}\langle X,X\rangle
    +\frac{a}{2}\langle X,X\rangle.
\end{eqnarray*}
This yields $(d/dt)\langle X,X\rangle=a\langle X,X\rangle$ and
since $\langle X(0),X(0)\rangle=0$ we obtain $\langle
X(t),X(t)\rangle=0$ for all $t$ and thus $X=0$.

Using the definition of $\gamma$ we obtain
\[
\frac{d}{dt}\langle\dot\gamma,\dot\gamma\rangle
=2\langle\nabla_{\dot\gamma}\dot\gamma,\dot\gamma\rangle
=a\langle\dot\gamma,\dot\gamma\rangle.
\]

The definition of $\gamma$, the fact that $B$ is geodesic and
\eqref{eqDXB} yield
\[
\frac{da}{dt}=\frac{d}{dt}\langle\dot\gamma,B\rangle
=\langle\nabla_{\dot\gamma}\dot\gamma,B\rangle
    +\langle\nabla_{\dot\gamma}B,\dot\gamma\rangle
=\frac{1}{2}\langle\dot\gamma,\dot\gamma\rangle
    -\frac{1}{2}\langle P\eta,\dot\gamma\rangle
=\frac{1}{2}\langle\dot\gamma,\dot\gamma\rangle
    -\frac{1}{2}\langle P\eta,P\eta\rangle.
\]

Now we calculate $(d/dt)\langle P\eta,P\eta\rangle$. Let
$\xi,\zeta\in\Gamma(\nu M)$ and $Y\in\Gamma(\g{C})$. Since $\g{C}$
is autoparallel, we have $\langle\enabla_BP\xi,Y\rangle=0$. Using
\eqref{II} we obtain
$\langle\enabla_BP\xi,\zeta\rangle=\langle\II(B,P\xi),\zeta\rangle=0$.
On the other hand we have $JP\xi=-\xi-PF\xi-F^2\xi
=-PF\xi-\sin^2(\varphi)\xi$, which gives, using~\eqref{II},
\begin{eqnarray*}
\langle\enabla_BP\xi,P\zeta\rangle
&=& \langle\enabla_BP\xi,J\zeta-F\zeta\rangle
    \ =\ -\langle\enabla_BJP\xi,\zeta\rangle
    -\langle\II(B,P\xi),F\zeta\rangle\\
&=& \langle\II(B,PF\xi),\zeta\rangle
    +\sin^2(\varphi)\langle\nabla_B^\perp\xi,\zeta\rangle
    \ =\ \langle P\nabla_B^\perp\xi,P\zeta\rangle.
\end{eqnarray*}
This readily implies,
\begin{equation}\label{eqDBPxi}
\enabla_BP\xi=P\nabla_B^\perp\xi \quad
\mbox{for all }\xi\in\Gamma(\nu M).
\end{equation}

Using again \eqref{II} we get
\begin{eqnarray*}
2\langle \enabla_{P\xi}P\xi ,Y\rangle
&=& -2\langle \enabla_{P\xi}Y, J\xi-F\xi \rangle
\,=\,2\langle JY,Z \rangle\langle \II(P\xi,Z),\xi \rangle
    + 2\langle Y,Z \rangle\langle \II(P\xi,Z),F\xi \rangle \\
&=& -\sin^2(\varphi)\langle JZ,Y \rangle\langle \xi,\xi \rangle
    +\sin^2(\varphi) \langle Y,Z \rangle \langle \xi,F\xi \rangle \
=\  \langle P\xi,P\xi \rangle\langle B,Y \rangle.
\end{eqnarray*}
Clearly, equation \eqref{II} implies
$\langle\enabla_{P\xi}P\xi,\zeta\rangle
=\langle\II(P\xi,P\xi),\zeta\rangle=0$. Using \eqref{II} and the
fact that $JP\xi=-PF\xi-\sin^2(\varphi)\xi$ we obtain
\begin{eqnarray*}
\langle\enabla_{P\xi}P\xi,P\zeta\rangle
&=& \langle\enabla_{P\xi}P\xi,J\zeta-F\zeta\rangle
    \ =\ -\langle\enabla_{P\xi}JP\xi,\zeta\rangle
    -\langle\II(P\xi,P\xi),F\zeta\rangle\\
&=& \langle\II(P\xi,PF\xi),\zeta\rangle
    +\sin^2(\varphi)\langle\enabla_{P\xi}\xi,\zeta\rangle
    \ =\ \langle P\nabla_{P\xi}^\perp\xi,P\zeta\rangle.
\end{eqnarray*}
Altogether this implies,
\begin{equation}\label{eqDPxiPxi}
\enabla_{P\xi}P\xi=\frac{1}{2}\langle P\xi,P\xi\rangle B
+P\nabla_{P\xi}^\perp\xi \quad
\mbox{for all }\xi\in\Gamma(\nu M).
\end{equation}

Finally, equations \eqref{eqDBPxi} and \eqref{eqDPxiPxi} yield
\begin{eqnarray*}
\frac{d}{dt}\langle P\eta,P\eta\rangle
&=& \frac{d}{dt}\langle\dot\gamma,P\eta\rangle
    \ =\ \langle\nabla_{\dot\gamma}\dot\gamma,P\eta\rangle
    +\langle\nabla_{\dot\gamma}P\eta,\dot\gamma\rangle\\
&=& \frac{a}{2}\langle P\eta,P\eta\rangle
    +a\langle P\nabla_B^\perp\eta,P\eta\rangle
    +\langle P\nabla_{P\eta}^\perp\eta,P\eta\rangle\\
&=& \frac{a}{2}\langle P\eta,P\eta\rangle
    +\langle P\nabla_{\dot\gamma}^\perp\eta,P\eta\rangle
\ =\ \frac{a}{2}\langle P\eta,P\eta\rangle
    +\frac{1}{2}\frac{d}{dt}\langle P\eta,P\eta\rangle.
\end{eqnarray*} and hence
\[
\frac{d}{dt}\langle P\eta,P\eta\rangle=a\langle P\eta,P\eta\rangle.
\]

Putting $b=\langle\dot\gamma,\dot\gamma\rangle$ and $c=\langle
P\eta,P\eta\rangle$ we then have the initial value problem:
\[
a'=\frac{1}{2}(b-c),\quad b'=ab,\quad c'=ac,\quad a(0)=0,\quad
b(0)=c(0)=1.
\]
Again, the uniqueness of solutions to differential equations yields
$a(t)=0$, $b(t)=c(t)=1$ for all $t$. Therefore,
$\langle\dot\gamma(t),\dot\gamma(t)\rangle=1$ and $\dot\gamma(t)\in
P\nu M$ for all $t$ as desired.

Let us assume then that $\gamma:I\to M$ is a curve satisfying
equation \eqref{eqGamma}. There exists a unit normal vector field
$\eta$ of $M$ in a neighborhood of $p$ such that
$\dot\gamma(t)=P\eta_{\gamma(t)}/\sin(\varphi)$ for all sufficiently
small $t$. Since $B$ is nonsingular and $\gamma$ is normal to $B$,
there exists a hypersurface $N$ in $M$ containing $\gamma$ and
transversal to $B$ in a neighborhood of $p$. The restriction of
$\eta$ to $N$ is a smooth unit normal vector field along $N$. We
define $\xi$ as the unit normal vector field on a neighborhood of
$p$ satisfying $\xi=\eta$ on $N$ and such that $\xi$ is obtained by
$\nabla^\perp$-parallel translation along the integral curves of
$B$. The smooth dependance on initial conditions of ordinary
differential equations implies that $\xi$ is smooth. Also, note that
$\nabla_B^\perp\xi=0$ and that $\xi$ is a local unit vector field
extending $\xi_p\in\nu_pM$.

The definition of $\xi$ and equations \eqref{eqDXB} and
\eqref{eqDBPxi} yield
$[B,P\xi]=\enabla_BP\xi-\enabla_{P\xi}B=\frac{1}{2}P\xi$, and hence
the distribution generated by $B$ and $P\xi$ is integrable. Let $U$
denote the integral submanifold through $p$. We prove that $U$ is an
open part of a totally geodesic $\R H^2\subset\CH$.

Since $B$ is geodesic we have $\enabla_BB=0$. Equation \eqref{eqDXB}
implies $2\enabla_{P\xi}B=-P\xi$, and by definition of $\xi$ we have
using \eqref{eqDBPxi} that $\enabla_BP\xi=P\enabla_B^\perp\xi=0$.
Now we prove that $2\enabla_{P\xi}P\xi=\langle P\xi,P\xi\rangle B$.

Let $\eta\in\nu M$ and denote by $\shape_\eta$ the shape operator of
$M$ with respect to the normal vector $\eta$. Equation \eqref{II}
implies that $\shape_\eta B=0$ for all $\eta$, and thus, for any
$\eta,\zeta\in\nu M$ the Ricci equation of $M$ reads
\[
\langle R^\perp_{B P\xi}\eta,\zeta\rangle
=\langle\eR(B,P\xi)\eta,\zeta\rangle
+\langle[\shape_\eta,\shape_\zeta]B,P\xi\rangle=0,
\]
where $R^\perp$ denotes the curvature tensor of the normal
connection $\nabla^\perp$. The previous equation, $2[B,P\xi]=P\xi$,
and the definition of $\xi$ imply
\[
0=R^\perp_{B P\xi}\xi
=\nabla_B^\perp\nabla_{P\xi}^\perp\xi
    -\nabla_{P\xi}^\perp\nabla_B^\perp\xi
    -\nabla_{[B,P\xi]}^\perp\xi
=\nabla_B^\perp\nabla_{P\xi}^\perp\xi
    -\frac{1}{2}\nabla_{P\xi}^\perp\xi,
\]
that is,
\begin{equation}\label{eqDperpPxixi}
2\nabla_B^\perp\nabla_{P\xi}^\perp\xi=\nabla_{P\xi}^\perp\xi.
\end{equation}

By definition of $\xi$, we have along $\gamma$ that
$2\enabla_{P\xi}P\xi=2\sin^2(\varphi)\enabla_{\dot\gamma}\dot\gamma
=2\sin^2(\varphi)\nabla_{\dot\gamma}\dot\gamma =\langle
P\xi,P\xi\rangle B$. On the other hand, \eqref{eqDPxiPxi} yields
$2\enabla_{P\xi}P\xi=\langle P\xi,P\xi\rangle
B+2P\nabla_{P\xi}^\perp\xi$, and hence $\nabla_{P\xi}^\perp\xi=0$
along $\gamma$. Now, let us take $\alpha$ an integral curve of $B$
through $\alpha(0)=\gamma(s)$. We have just seen that
$\nabla_{P\xi}^\perp\xi\,_{\vert_{\alpha(0)}}
=\nabla_{P\xi}^\perp\xi \,_{\vert_{\gamma(s)}}=0$. Moreover, using
\eqref{eqDperpPxixi} and the fact that $\shape_\eta B=0$ for any
$\eta\in\nu M$, we obtain
\[
2\enabla_{\dot\alpha}\nabla_{P\xi}^\perp\xi\,_{\vert_{t}}
=2\nabla_B^\perp\nabla_{P\xi}^\perp\xi\,_{\vert_{\alpha(t)}}
    -2\shape_{\nabla_{P\xi}^\perp\xi}B\,_{\vert_{\alpha(t)}}
=\nabla_{P\xi}^\perp\xi\,_{\vert_{\alpha(t)}}.
\]
Therefore, by the uniqueness of solutions to differential equations
we get $\nabla_{P\xi}^\perp\xi\,_{\vert_{\alpha(t)}}=0$ for all $t$,
and as a consequence $2\enabla_{P\xi}P\xi=\langle P\xi,P\xi\rangle
B$ along the integral submanifold $U$. Hence, $U$ is an open part of
a totally geodesic $\R H^2\subset\CH$.

We define $\bar P\xi = P\xi/\lVert P\xi \rVert = P\xi/\sin(\varphi)$
along $U$. From \eqref{eqDPxiPxi} we obtain $2\enabla_{\bar
P\xi}\bar P\xi = B$ since $\xi$ is unit normal. Using this and
\eqref{eqDXB} we get
\[
\enabla_{\bar P\xi}\enabla_{\bar P\xi}\bar P\xi+ \langle
\enabla_{\bar P\xi}\bar P\xi, \enabla_{\bar P\xi}\bar P\xi\rangle
\bar P\xi= \frac{1}{2}\enabla_{\bar P\xi}B + \frac{1}{4}\langle
B,B\rangle\bar P\xi= 0.
\]
From this we see that the integral curves of $\bar P\xi$ are
horocycles with center $x$ at infinity contained in an open part of
a totally geodesic real hyperbolic plane contained in $\CH$.
Corollary \ref{thConstruction}, the rigidity of totally geodesic
submanifolds of Riemannian manifolds (see e.g.\ \cite{BCO03}, p.\
230), and of horocycles in real hyperbolic planes (see e.g.\
\cite{BCO03}, pp.\ 24-26), then imply the assertion.
\end{proof}

\begin{remark}
{\rm The proof shows that the differential equation (\ref{diffeq})
characterizes left translation of the normal spaces by
$S_{\g{c}}$.}
\end{remark}

\section{The tubes around $W^{2n-k}$ and $W^{2n-k}_{\varphi}$}

To accomplish the task of investigating the geometry of orbits of
the cohomogeneity one actions on $\CH$ we will deal with two
different possibilities depending on the constant K\"{a}hler angle
$\varphi\in(0,\pi/2)$ or $\varphi = \pi/2$ of $\g{w}$. For this we
first we recall a few properties of the solvable foliation already
studied in \cite{Be98}.

\subsection{The solvable foliation}\label{scSolvable}

The solvable foliation is the foliation on $\CH$ arising from $k=1$.
In this case $\varphi=\pi/2$ and the orbit $S\cdot o=\W[]{2n-1}$ is
a minimal homogeneous ruled real hypersurface. Its principal
curvatures are $1/2$, $-1/2$ and $0$ with multiplicities $1$, $1$
and $2n-3$. The following theorem shows that this eigenvalue
structure is characteristic of this orbit.

\begin{theorem}[Rigidity of the submanifold $W^{2n-1}$]
Let $M$ be a
connected real hypersurface in $\CH$, $n \geq 2$,
with three distinct principal curvatures $1/2$, $-1/2$ and $0$ and
multiplicities $1$, $1$ and $2n-3$, respectively. Then $M$ is
holomorphically congruent to an open part of the minimal homogeneous
ruled real hypersurface $\W[]{2n-1}$.
\end{theorem}

This result was proved in \cite{BD05} for $n\geq 3$. The analogous
statement for $n=2$ is more involved and follows from the
classification of real hypersurfaces with constant principal
curvatures in the complex hyperbolic plane \cite{BD06}. Any other
orbit of the action of $S$ is an equidistant hypersurface to this
minimal one. Any two such orbits are congruent to each other if and
only if their distance to $S\cdot o = W^{2n-1}$ is the same. None of
them is ruled by a totally geodesic $\CH[n-1]$. Let $M(r)$ denote an
orbit of $S$ at a distance $r>0$ from $S\cdot o$. The shape operator
of $M(r)$ has exactly three eigenvalues
\[
\lambda_{1/2}=\frac{3}{4}\tanh\left(\frac{r}{2}\right)
\pm \frac{1}{2}\sqrt{1-\frac{3}{4}\tanh^2\left(\frac{r}{2}\right)}\ ,
\quad
\lambda_3=\frac{1}{2}\tanh\left(\frac{r}{2}\right)\ .
\]
with corresponding multiplicities $m_1=1$, $m_2=1$ and $m_3=2n-3$.
The Hopf vector field $J\xi$ has nontrivial projection onto the
principal curvature spaces of $\lambda_1$ and $\lambda_2$.

The subspace $\g{a} + \g{w}^\perp + J\g{w}^\perp + \g{g}_{2\alpha}$
of $\g{a} + \g{n}$ is a subalgebra of $\g{a}+\g{n}$, and the orbit
through $o$ of the corresponding connected closed subgroup of $AN$
is a totally geodesic  $\CH[2]$. The action of the connected closed
subgroup of $S$ with Lie algebra $\g{a} + J \g{w}^\perp +
\g{g}_{2\alpha}$ induces the solvable foliation on this totally
geodesic $\CH[2]$. The relevant geometric information of the
solvable foliation on $\CH$ is contained in the ``slice'' $\CH[2]$.
We describe in more detail the geometry of the leaves of the
solvable foliation on $\CH[2]$ in what follows.

Let $\gamma$ be the geodesic in $\CH[2]$ determined by $\gamma(0) =
o$ and $\dot\gamma(0)=\xi$, where $\xi \in \g{w}^\perp$ is a unit
vector. Let $r\in\R$ and denote by $M(r)$ the leaf of the solvable
foliation containing $\gamma(r)$. According to (\ref{geodesictangentvector})
the tangent vector field
$\dot{\gamma}$ of the geodesic $\gamma$ can be written with respect
to left-invariant vector fields as
\[
\dot\gamma(t) =-\tanh(t/2)B+\sech(t/2)\xi.
\]
Then $\{Z,J\xi,\sech(r/2)B+\tanh(r/2)\xi\}$ is an orthonormal
basis of $T_{\gamma(r)}M(r)$. The distribution on
$M(r)$ generated by $Z$ and $J\xi$ is integrable
by \eqref{bracket}. Moreover, using
\eqref{LeviCivita} we get
\begin{equation}\label{eqNablaZJxi}
\enabla_Z Z=B\ ,\ 2\enabla_{J\xi}J\xi=B\ ,
\ 2\enabla_ZJ\xi=\xi\ ,\ 2\enabla_{J\xi}Z=\xi.
\end{equation}
Thus the shape operator of the leaf of this distribution containing
$\gamma(r)$ with respect to the unit normal vector
$\sech(r/2)B_{\gamma(r)}+\tanh(r/2)\xi_{\gamma(r)}\in
T_{\gamma(r)}M(r)$ is given by the matrix
\[
\frac{1}{2}\left(
\begin{array}{cc}
2\sech\left(\frac{r}{2}\right) & \tanh\left(\frac{r}{2}\right)\\[1ex]
\tanh\left(\frac{r}{2}\right)  & \sech\left(\frac{r}{2}\right)
\end{array}\right)
\]
with respect to the basis $\{Z_{\gamma(r)},J\xi_{\gamma(r)}\}$.
Using the Gauss equation and \eqref{eqNablaZJxi} we get that the Gaussian
curvature of the leaf through $\gamma(r)$ is equal to zero. For
topological reasons it is clear that the leaf is a Euclidean plane
$\R^2$.

On the other hand, using Lemma \eqref{LeviCivita} we get
\[
\enabla_{\sech\left(\frac{r}{2}\right)B
+\tanh\left(\frac{r}{2}\right)\xi}
\left(\sech\left(\frac{r}{2}\right)B
+\tanh\left(\frac{r}{2}\right)\xi\right)
=-\frac{1}{2}\tanh\left(\frac{r}{2}\right)
\dot\gamma(r).
\]
Hence, every integral curve of $\sech(r/2)B+\tanh(r/2)\xi$ is a
geodesic in $M(r)$.

All in all, this means

\begin{theorem}
The leaves of the solvable foliation on $\CH[2]$ are diffeomorphic
to $\R^3$ and are foliated orthogonally by a one-dimensional totally
geodesic foliation and a two-dimensional foliation whose leaves are
Euclidean planes.
\end{theorem}

\subsection{Constant K\"{a}hler angle
$\varphi=\pi/2$}\label{scKahlerPi2}

In this case $\g{w}^\perp$ has constant K\"{a}hler angle
$\varphi=\pi/2$, that is, $\g{w}^\perp$ is real. This means that the
normal bundle $\nu\W[]{2n-k}$ of $\W[]{2n-k}$ is totally real. We
recall that the second fundamental form of $\W[]{2n-k}$ is given by
the trivial symmetric bilinear extension of $\II(Z,J\xi)=(1/2)\xi$
for all $\xi\in\nu\W[]{2n-k}$. Thus the eigenvalues of the shape
operator of $\W[]{2n-k}$ with respect to any unit vector
$\xi\in\nu\W[]{2n-k}$ are $1/2$, $-1/2$ and $0$ with multiplicities
$1$, $1$ and $2n-2-k$ respectively. The corresponding principal
curvature spaces are $\R(Z+J\xi)$, $\R(Z-J\xi)$ and
$T\W[]{2n-k}\ominus(\R Z + \R J\xi)$ respectively.

The above information allows us to calculate the shape operator of
the principal orbits using Jacobi field theory. Every principal
orbit of this action is a tube around $\W[]{2n-k}$. We denote by
$M(r)$ the tube at distance $r > 0$ and fix $o\in\W[]{2n-k}$ and a
unit vector $\xi\in\nu_o\W[]{2n-k}$. Let $\gamma_\xi$ be the
geodesic in $\CH$ given by the initial conditions $\gamma_\xi(0)=o$
and $\dot{\gamma}_\xi(0)=\xi$. We recall that the Jacobi equation in
the complex hyperbolic space of constant holomorphic sectional
curvature $-1$ along $\gamma_\xi$ reads $4\zeta_X''-\zeta_X-3\langle
J\dot\gamma_\xi,\zeta_X\rangle J\dot\gamma_\xi=0$.

For any $X\in T_o\CH$ we denote by $\pd_X$ the parallel displacement
of the vector $X$ along $\gamma_\xi$. If $X\in T_o\W[]{2n-k}$ we
denote by $\zeta_X$ the Jacobi field along $\gamma_\xi$ defined by
the initial conditions $\zeta_X(0)=X$ and
$\zeta_X'(0)=-\shape_\xi(X)$. If $X$ is a principal curvature
vector, that is, $\shape_\xi X=\lambda X$ for some $\lambda\in\R$,
then the Jacobi equation can be solved explicitly to get
\[
\zeta_X(t)=f_\lambda(t)\pd_X(t)+\langle X,J\xi\rangle g_\lambda(t)
J\dot\gamma_\xi(t)
\]
with
\[
f_\lambda(t)=\cosh\Bigl(\frac{t}{2}\Bigr)-2\lambda\sinh
    \Bigl(\frac{t}{2}\Bigr),\
g_\lambda(t)=\left(\cosh\Bigl(\frac{t}{2}\Bigr)-1\right)
    \left(1+2\cosh\Bigl(\frac{t}{2}\Bigr)
    -2\lambda\sinh\Bigl(\frac{t}{2}\Bigr)\right).
\]
If $X\in\nu_o\W[]{2n-k}\ominus\R\xi$ we define the Jacobi field
$\zeta_X$ along $\gamma_\xi$ by the initial conditions
$\zeta_X(0)=0$ and $\zeta_X'(0)=X$. In this case we have
\[
\zeta_X(t)=p(t)\pd_X(t)+\langle X,J\xi\rangle q(t) J\dot\gamma_\xi(t)
\]
with
\[
p(t)=2\sinh\Bigl(\frac{t}{2}\Bigr)\ ,\
q(t)=2\sinh\Bigl(\frac{t}{2}\Bigr)
    \left(\cosh\Bigl(\frac{t}{2}\Bigr)-1\right).
\]

Using the above formulas one easily gets
\[
\zeta_X(t)=\left\{ {\everymath{\displaystyle}
\begin{array}{l@{\quad,\qquad\mbox{if }}l}
\cosh\left(\frac{t}{2}\right)\pd_Z(t)
-\frac{1}{2}\sinh(t)\pd_{J\xi}(t)
    & X=Z,\\[2ex]
-\sinh\left(\frac{t}{2}\right)\pd_Z(t)+\cosh(t)\pd_{J\xi}(t)
    & X=J\xi,\\[1ex]
\cosh\left(\frac{t}{2}\right)\pd_X(t)
    & X\in T\W[]{2n-k}\ominus(\R Z+\R J\xi),\\[1ex]
2\sinh\left(\frac{t}{2}\right)\pd_X(t)
    & X\in\nu\W[]{2n-k}\ominus\R\xi.
\end{array}}
\right.
\]

We define the endomorphism $D(r)$ of
$T_{\gamma_\xi(r)}M(r)\ominus\R\dot\gamma_\xi(r)$ by
$D(r)B_X(r)=\zeta_X(r)$ for all $X\in T_{o}\CH\ominus\R\xi$. Jacobi
field theory shows that the shape operator of $M(r)$ at
$\gamma_\xi(r)$ with respect to  $-\dot\gamma(r)$ is given by
$\shape(r)=D'(r)D(r)^{-1}$. In our case $\shape(r)$ is represented
by the matrix
\[
\renewcommand{\arraystretch}{1.25}
\shape(r)=\frac{1}{2}\left(
\begin{array}{@{}cc|c|c@{}}
\tanh^3\left(\frac{r}{2}\right)
    &   -\sech^3\left(\frac{r}{2}\right) &&\\
-\sech^3\left(\frac{r}{2}\right)
    &   2\left(1+\frac{1}{2}\sech^2\left(\frac{r}{2}\right)\right)
        \tanh\left(\frac{r}{2}\right) &&\\
\hline
&&  \tanh\left(\frac{r}{2}\right)\,\Id_{2n-2-k}  &\\
\hline
&&& \coth\left(\frac{r}{2}\right)\,\Id_{k-1}
\end{array}\right)\!.
\]
with respect to the orthogonal sum decomposition
\[
T_{\gamma_\xi(r)}M(r)=\pd_{\R Z+\R J\xi}(r) +
\pd_{T_o\W[]{2n-k}\ominus(\R Z+\R J\xi)}(r)+
\pd_{\nu_o\W[]{2n-k}\ominus\R\xi}(r),
\]
where $\pd_V$ denotes the parallel translation of any vector
subspace $V\subset T_o\CH$ along $\gamma_\xi$.

A straightforward calculation shows that $M(r)$ has four principal
curvatures
\[
{\everymath{\displaystyle}
\begin{array}{@{}rcl@{\quad}rcl@{}}
\lambda_1&=&\frac{3}{4}\tanh\left(\frac{r}{2}\right)
-\frac{1}{2}\sqrt{1-\frac{3}{4}\tanh^2\left(\frac{r}{2}\right)}\ ,&
\lambda_2&=&\frac{3}{4}\tanh\left(\frac{r}{2}\right)+
\frac{1}{2}\sqrt{1-\frac{3}{4}
\tanh^2\left(\frac{r}{2}\right)}\ ,\\[2ex]
\lambda_3&=&\frac{1}{2}\tanh\left(\frac{r}{2}\right)\ ,&
\lambda_4&=&\frac{1}{2}\coth\left(\frac{r}{2}\right)\
\end{array}}
\]
with corresponding multiplicities $m_1 = m_2=1$, $m_3=2n-2-k$ and
$m_4=k-1$. The Hopf vector field on $M$ has nontrivial orthogonal
projection onto the principal curvature spaces of $\lambda_1$ and
$\lambda_2$. A special situation occurs when $r=\ln(2+\sqrt{3})$. In
this case we have $\lambda_2=\lambda_4$ and the principal curvatures
are $\lambda_1=0$, $\lambda_2=\lambda_4=\sqrt{3}/2$ and
$\lambda_3=\sqrt{3}/6$ with multiplicities $1$, $k$ and $2n-k-2$
respectively.

The previous calculations show that the interesting part of the
shape operator of both the singular orbit $\W[]{2n-k}$ and the
principal orbit $M(r)$ concerns the vectors $Z$ and $J\xi$. More
precisely, let $\xi\in \nu_o\W[]{2n-k}$ be a unit vector. Consider
the subalgebra $\tilde{\g{g}}=\g{a} + \R\xi + \R J\xi +
\g{g}_{2\alpha}$ of $\g{a}+\g{g}_\alpha+\g{g}_{2\alpha}$, and let
$\tilde{G}$ be the connected closed subgroup of $A N$ with Lie
algebra $\tilde{\g{g}}$. The orbit $\tilde{G}\cdot o$ is a totally
geodesic $\CH[2]$ in $\CH$. This $\CH[2]$ defines a ``slice'' of the
action of $N^0_K(S)S$ through $o$. Next,
$\tilde{\g{h}}=\g{s}\cap\tilde{\g{g}}$ is a subalgebra of
$\tilde{\g{g}}$ of codimension one. Let $\tilde{H}$ be the connected
closed subgroup of $\tilde{G}$ with Lie algebra $\tilde{\g{h}}$.
Then $\tilde{H}$ acts on $\CH[2]=\tilde{G}\cdot o$ with
cohomogeneity one and gives exactly the solvable foliation of
$\CH[2]$ described in Section \ref{scSolvable}. The orbits of the
action of $\tilde{H}$ on $\CH[2]$ are the equidistant hypersurfaces
to the orbit $\tilde{H}\cdot o$. On the other hand, the intersection
of the orbits of the cohomogeneity one action of $N^0_K(S)S$ on
$\CH$ with the slice $\CH[2]$ also gives tubes around
$\tilde{H}\cdot o$ because $\CH[2]=\tilde{G}\cdot o$ is totally
geodesic in $\CH$. Thus, the geometry of the orbits of the action of
$G$ on $\CH$ in the slice $\CH[2]$ is exactly the geometry of the
orbits of the action of $\tilde{H}$ on $\CH[2]$. This study was
accomplished in the previous subsection.

\subsection{Constant K\"{a}hler angle
$\varphi\in(0,\pi/2)$}\label{scKahlerPhi}

Again, we assume the notation above, and consider the singular orbit
$\W{2n-k}$ of the cohomogeneity one action determined by the Lie
group $N_K^0(S)S$, where $S$ is the connected, simply connected Lie
group whose Lie algebra is $\g{s}=\g{a}+\g{w}+\g{g}_{2\alpha}$,
where $\g{w}^\perp$ has constant K\"ahler angle
$\varphi\in(0,\pi/2)$. In this case we have that $k$ is an even
number. The second fundamental form of $\W{2n-k}$ is given by the
trivial bilinear extension of $\II(Z,P\xi)=(\sin^2(\varphi)/2)\xi,$
for each unit $\xi\in\nu\W{2n-k}$. Thus, the eigenvalues of the
shape operator with respect to $\xi$  are $\sin(\varphi)/2$,
$-\sin(\varphi)/2$ and $0$, with multiplicities $1$, $1$ and
$2n-k-2$.

It is convenient to introduce the notation
\[
\bar{P}\xi=P\xi/\lVert P\xi\rVert=P\xi/\sin(\varphi)
\qquad\mbox{and}\qquad
\bar{F}\xi=F\xi/\lVert F\xi\rVert=F\xi/\cos(\varphi)
\]
for each unit vector $\xi\in\nu\W{2n-k}$. Then, the eigenspaces of
$\sin(\varphi)/2$, $-\sin(\varphi)/2$ and $0$ of the shape operator
of $\W{2n-k}$ with respect to $\xi$ are $\R(Z+\bar{P}\xi)$,
$\R(-Z+\bar{P}\xi)$ and $T\W{2(n-k)}\ominus(\R Z+\R\bar{P}\xi)$,
respectively.

The shape operator of the principal orbits can be calculated using
Jacobi field theory as in the previous section. We delete the
calculations, which are straightforward (although long) and directly
give the matrix representation of the shape operator $\shape(r)$
in direction ${-\dot{\gamma}_\xi(r)}$ of the orbit at
a distance $r>0$ from $\W{2n-k}$
\[
\shape(r) = \left(
{\renewcommand{\arraystretch}{1.5}
\begin{array}{c|c|c}
s(r)&&\\
\hline
& \frac{1}{2}\tanh\left(\frac{r}{2}\right)\,\Id_{2n-k-2}&\\
\hline &&\frac{1}{2}\coth\left(\frac{r}{2}\right)\;\Id_{k-2}
\end{array}}\right),
\]
with respect to the direct sum decomposition
\[
T_{\gamma_\xi(r)}M(r)=
    \pd_{\R Z+\R\bar{P}\xi+\R\bar{F}\xi}(r)
    +\pd_{T\W{2n-k}\ominus(\R Z+\R\bar{P}\xi)}(r)
    +\pd_{\nu\W{2n-k}\ominus(\R\xi+\R\bar{F}\xi)}(r).
\]
Here, $s(r)$ is a symmetric $3\times 3$ real matrix whose explicit
entries we do not provide (they can be obtained after some
elementary but long calculations). The characteristic polynomial of
$s(r)$ is
\[
\everymath{\displaystyle}
\begin{array}{@{}r@{\ }c@{\ }l@{}}
p_{r,\varphi}(x)
&=&-x^3+\frac{1}{2}\left\{\csch\Bigl(\frac{r}{2}\Bigr)
\sech\Bigl(\frac{r}{2}\Bigr)
\!+4\tanh\Bigl(\frac{r}{2}\Bigr)\right\}x^2
-\frac{1}{4}\left\{2\sech^2\Bigl(\frac{r}{2}\Bigr)
+5\tanh^2\Bigl(\frac{r}{2}\Bigr)\right\}x\\
&&-\frac{1}{8}\csch\Bigl(\frac{r}{2}\Bigr)
\sech^3\Bigl(\frac{r}{2}\Bigr)\left\{
\sin^2(\varphi)-\sinh^2\Bigl(\frac{r}{2}\Bigr)
-2\sinh^4\Bigl(\frac{r}{2}\Bigr)\right\}.
\end{array}
\]
If we introduce the variable
$6x=\coth(r/2)\,z-\csch(r/2)\sech(r/2)-4\tanh(r/2)$, then the
polynomial equation $p_{r,\varphi}(x)=0$ transforms into
$z^3-3z+\beta_{r,\varphi}=0$, where
$\beta_{r,\varphi}=27\sin^2(\varphi)\tanh^2(r/2)\sech^4(r/2)-2$. The
discriminant of this cubic equation is
$\Delta_{r,\varphi}=27(\beta_{r,\varphi}^2-4)$. It is easy to prove
that $\Delta_{r,\varphi}<0$ for all $r>0$, which means that the
above cubic equation has exactly three distinct real roots for any
$r$. They can be calculated explicitly as follows. Let
$u_{r,\varphi}^i$, $i\in\{1,2,3\}$, denote each cubic root of the
unit complex number
$(\beta_{r,\varphi}+\sqrt{\beta_{r,\varphi}^2-4})/2$. Then,
$-u_{r,\varphi}^i-1/u_{r,\varphi}^i$ is a solution to
$z^3-3z+\beta_{r,\varphi}=0$ and hence, the eigenvalues of $s(r)$
are given by
\[
\lambda_i(r)=-\frac{1}{6} \left(\coth\Bigl(\frac{r}{2}\Bigr)
    \Bigl(u_{r,\varphi}^i+\frac{1}{u_{r,\varphi}^i}\Bigr)
+\csch\Bigl(\frac{r}{2}\Bigr)\sech\Bigl(\frac{r}{2}\Bigr)
+4\tanh\Bigl(\frac{r}{2}\Bigr)\right),\ i\in\{1,2,3\}.
\]
On the other hand, $p_{r,\varphi}((1/2)\tanh(r/2))\neq 0$ and
$p_{r,\varphi}((1/2)\coth(r/2))\neq 0$. Thus, neither
$(1/2)\tanh(r/2)$ nor $(1/2)\coth(r/2)$ are eigenvalues of $s(r)$.
This implies that $M(r)$ has five distinct constant principal
curvatures when $k>2$ and four distinct principal curvatures when
$k=2$.

Hereafter, we follow the procedure of the previous section and focus
our study on the non-trivial part of the shape operator of $M(r)$.
Let $\g{v}_0\subset\g{g}_\alpha$ be a two-dimensional vector
subspace with constant K\"{a}hler angle $\varphi$. Then,
$\tilde{\g{g}}=\g{a}+\C\,\g{v}_0+\g{g}_{2\alpha}$ is a Lie
subalgebra of $\g{a}+\g{n}$. Let $\tilde{G}=\Exp(\tilde{\g{g}})$ be
the connected, simply connected Lie subgroup of $A N$ whose Lie
algebra is $\tilde{\g{g}}$. Then, $\tilde{G}\cdot o$ is a totally
geodesic $\CH[3]$ in $\CH$ containing $o$. The vector subspace
$\tilde{\g{h}}=\g{a}+\g{v}_0+\g{g}_{2\alpha}$ is a Lie subalgebra of
$\tilde{\g{g}}$ of codimension two. Denote by
$\tilde{H}=\Exp(\tilde{\g{h}})$ the connected, simply connected Lie
subgroup of $\tilde{G}$ whose Lie algebra is $\tilde{\g{h}}$. We
know that the Lie group $N_K^0(\tilde{H})\tilde{H}$ acts on
$\tilde{G}\cdot o$ with cohomogeneity one and its orbit through $o$
is exactly $\tilde{H}\cdot o$. This cohomogeneity one action is the
one we have been describing throughout this subsection. We are
interested in this particular case because it is the simplest of all
cases containing all of the interesting geometry of the tubes.

Let $M(r)$ denote the tube around $\W{4}\subset\CH[3]$ at distance
$r>0$. Then, $M(r)$ is the principal orbit of the action
$N_K^0(\tilde{H})\tilde{H}$ at a distance $r$ from the singular
orbit $\tilde{H}\cdot o=\W{4}$. The normal exponential map
$\exp^\perp:\nu\W{4}\to\CH[3]$ of $\W{4}$ is a diffeomorphism and
hence for each $p\in M(r)$ there exists a unique unit vector
$\xi(p)\in\nu\W{4}$ such that $p=\exp(r\,\xi(p))$. Clearly, the map
$p\mapsto\xi(p)$ is differentiable. Let
$\gamma_{\xi(p)}(t)=\exp(t\,\xi(p))$ be the unique geodesic
perpendicular to $\W{4}$ that joins $\W{4}$ and $p$. For any $X\in
T_{\gamma_{\xi(p)}(0)}\CH[3]$ we denote by $\pd_X^p(r)$ the parallel
displacement of $X$ to the point $p$ along the geodesic
$\gamma_{\xi(p)}$. The smooth dependence on initial conditions of
the solution to ordinary differential equations implies that
$\pd_X(r):p\mapsto\pd_X^p(r)$ is a smooth vector field on $\CH[3]$.
Moreover, if $X$ is tangent to $\W{4}$, then $\pd_X^p(r)$ is tangent
to $M(r)$. We have

\begin{theorem} The following two statements hold.
\begin{itemize}
\item[(i)] Let ${\mathcal D}$ be the rank one distribution on $M(r)$
defined by $\pd_B^p(r)$, $p \in M(r)$, and denote by ${\mathcal
D}^\perp$ the orthogonal complement of ${\mathcal D}$ in $TM(r)$.
Then both ${\mathcal D}$ and ${\mathcal D}^\perp$ are integrable.
Moreover, ${\mathcal D}$ is autoparallel, that is, each of its
leaves is totally geodesic in $M(r)$. If $p \in M(r)$ and $\R H^2$
is the totally geodesic real hyperbolic plane in $\C H^3$ which is
determined by $\xi(p)$ and $B_o$, where $o \in \W{4}$ is the
footpoint of $\xi(p)$, then the leaf of ${\mathcal D}$ through $p$
is parametrized by the parallel curve through $p$ in $\R H^2$ of the
geodesic in $\R H^2$ through $o$ and in direction $B_o$.

\item[(ii)] Let ${\mathcal E}$ be the rank two distribution on
$M(r)$ defined by $\R\pd_B^p(r) + \R\pd_{PF\xi(p)}^p(r)$. Then
${\mathcal E}$ is autoparallel and each integral manifold has
constant sectional curvature $-(1/4)\sech(r/2)$. If $p \in M(r)$ and
$\R H^3$ is the totally geodesic real hyperbolic 3-space in $\C H^3$
which is determined by $\xi(p)$, $PF\xi(p)$ and $B_o$, where $o \in
\W{4}$ is the footpoint of $\xi(p)$, then the leaf of ${\mathcal E}$
through $p$ is the parallel surface through $p$ in $\R H^3$ of the
totally geodesic $\R H^2$ in $\R H^3$ through $o$ determined by
$PF\xi(p)$ and $B_o$.
\end{itemize}
\end{theorem}

\begin{proof}
Let $p\in M(r)$ and $o\in\W{4}$ the footpoint of $\xi(p)$. The
vectors $B_o$ and $\xi(p)$ determine a totally geodesic real
hyperbolic plane $\R H^2\subset\CH[3]$ through $o$. Let $\tilde p\in
M(r)\cap\R H^2$. Since the normal exponential map of
$\W{4}\subset\CH[3]$ is a diffeomorphism and $\W{4}\cap\R H^2$ is
the path of the geodesic determined by $B_o$ we have that
$\xi(\tilde p)\in T_{\tilde o}\R H^2$, where $\tilde o\in\W{4}\cap\R
H^2$ is the footpoint of $\xi(\tilde p)$. Since $\R H^2$ is totally
geodesic, $\pd_B^{\tilde p}(r)$ is tangent to $M(r)\cap\R H^2$. This
proves that $M(r)\cap\R H^2$ is an integral manifold of ${\mathcal
D}$ through $p$. Moreover, if $X\in\Gamma(T(M(r)\cap\R H^2))$ and
$\nabla^{M(r)}$ denotes the Levi-Civita connection of $M(r)$, it is
clear that $\nabla^{M(r)}_XX\in\Gamma(TM(r))$. On the other hand,
since $\R H^2$ is totally geodesic, $\enabla_XX\in\Gamma(T\R H^2)$
and hence $\nabla^{M(r)}_XX\in\Gamma({\mathcal D})$, which proves
that ${\mathcal D}$ is autoparallel and the first part of (i)
follows.

Similarly, let $\R H^3\subset\CH[3]$ be the totally geodesic real
hyperbolic space determined by $B_o$, $\xi(p)$ and $PF\xi(p)$. Then
the integral submanifold of ${\mathcal E}$ through $p$ is
$M(r)\cap\R H^3$, and since $\R H^3$ is totally geodesic
and intersects $M(r)$ perpendicularly, we see
that ${\mathcal E}$ is autoparallel. The curvature of the integral
submanifolds of ${\mathcal E}$ is $-(1/4)\sech(r/2)$ as they are
equidistant to a totally geodesic $\R H^2\subset\R H^3$ obtained as
the intersection of $\W{4}$ and $\R H^3$. This proves (ii).

Now we prove that ${\mathcal D}^\perp$ is integrable. We define the
vector field $\tilde\xi$ along $\CH[3]\setminus\W{4}$ by
\[
\tilde\xi_{\exp(r\eta)}=L_{\exp(r\eta)*}\eta\quad\mbox{for all unit
vectors $\eta\in\nu\W{4}$ and $r>0$.}
\]
Let $\eta\in\nu\W{4}$ be a unit vector and denote also by $\eta$ the
unit vector field on $\CH[3]$ obtained by left translation to all
points of $\CH[3]$. Let $\gamma_\eta$ be the geodesic in $\CH[3]$ such that
$\dot\gamma_\eta(0)=\eta$. According to (\ref{geodesictangentvector}) we have
$\dot{\gamma}_\eta(r)=-\tanh(r/2)B_{{\gamma}_\eta(r)}
+\sech(r/2)\eta_{{\gamma}_\eta(r)}$ where $B$ and $\eta$ are
considered as left-invariant vector fields. Using \eqref{LeviCivita}
we get
\begin{eqnarray*}
\enabla_{\dot\gamma_\eta}\left(\sech\Bigl(\frac{r}{2}\Bigr)B
+\tanh\Bigl(\frac{r}{2}\Bigr)\eta\right) &=&
-\frac{1}{2}\sech\Bigl(\frac{r}{2}\Bigr)
    \tanh\Bigl(\frac{r}{2}\Bigr)B
    +\sech\Bigl(\frac{r}{2}\Bigr)\enabla_{\dot\gamma_\eta}B\\
&&  +\frac{1}{2}\sech^2\Bigl(\frac{r}{2}\Bigr)\eta
    +\tanh\Bigl(\frac{r}{2}\Bigr)\enabla_{\dot\gamma_\eta}\eta
\ =\ 0.
\end{eqnarray*}
This proves that
\begin{equation}\label{eqBB}
\pd_B^p(r)=\sech\Bigl(\frac{r}{2}\Bigr)B_p
    +\tanh\Bigl(\frac{r}{2}\Bigr)\tilde\xi_p.
\end{equation}

Now, let $p\in M(r)$. Let us assume without loss of generality that
$\gamma_{\xi(p)}(0)=o$ and write $\eta=\xi(p)\in\nu_oW^4_\varphi$. The
formulas for $\dot\gamma_\eta(r)$ and $\pd^p_B(r)$ show that
${\mathcal D}^\perp_p$ is spanned by $Z_p$, $P\eta_p$, $PF\eta_p$
and $F\eta_p$. Let $X,Y\in\Gamma({\mathcal D}^\perp)$. Using
\eqref{eqBB} we get
\[
\langle\enabla_{X_p}Y,\pd^p_B(r)\rangle
=-\langle Y_p,\enabla_{X_p}\pd_B(r)\rangle
=-\sech\Bigl(\frac{r}{2}\Bigr)\langle Y_p,\enabla_{X_p}B\rangle
-\tanh\Bigl(\frac{r}{2}\Bigr)
    \langle Y_p,\enabla_{X_p}\tilde\xi\rangle.
\]
The first term on the right-hand side of this equation may be
calculated using \eqref{LeviCivita} so we turn our attention to
$\enabla_{X_p}\tilde\xi$. Let $\chi_X$ be a curve such that
$\dot\chi_X(0)=X_p$. The curve $\chi_X$ can be written as
$\chi_X(t)=\exp_{g_X(t)\cdot
o}\bigl(s_X(t)(j_X(t)\bar{F}\eta+h_X(t)\eta)\bigr)$ for certain
smooth functions $g_X:I\to\tilde{H}$ and $s_X,j_X,h_X:I\to\R$
satisfying $s_X(0)=r$, $j_X(0)=0$, $h_X(0)=1$ and $j_X^2+h_X^2=1$.
Taking derivatives on the last equality we get $h_X'(0)=0$. Since
$\tilde\xi_{\chi_X(t)}=j_X(t)\bar{F}\eta_{\chi_X(t)}
+h_X(t)\eta_{\chi_X(t)}$, \eqref{LeviCivita} yields
\[
\enabla_{\dot\chi_X(0)}\tilde\xi
=\enabla_{\dot\chi_X(0)}\Bigl(j_X\bar{F}\eta+h_X\eta\Bigr)
=j_X'(0)\bar{F}\eta_p+\enabla_{\dot\chi_X(0)}\eta.
\]
Again, the second term can be calculated using \eqref{LeviCivita}.
All in all this means (interchanging the roles of $X$ and $Y$) that
\[
\langle[X,Y]_p,\pd^p_B(r)\rangle
=\langle\enabla_{X_p}Y-\enabla_{Y_p}X,\pd^p_B(r)\rangle
=-\tanh\Bigl(\frac{r}{2}\Bigr)\langle j_X'(0)Y_p-j_Y'(0)X_p,
    \bar{F}\eta_p\rangle.
\]
Note that the vector fields $Z,P\tilde{\xi},PF\tilde{\xi},F\tilde{\xi}$
restricted to $M(r)$ form a global frame field of ${\mathcal D}^\perp$,
and at the point $p$ we have $\tilde{\xi}_p = \eta_p$. Therefore
it is clear that the result follows if we
prove $j_X'(0)=0$ for $X_p\in\{Z_p,P\eta_p,PF\eta_p\}$.

The curve $\alpha(t)=\Exp_{\g{a}+\g{n}}(t X_o)$ is tangent to $X_o$
for $t=0$. Then, $\chi_X(t)=\gamma_\eta(r)\alpha(t)$ is tangent to
$X_p$ at $t=0$. We define
\[
U_Z=\sech^2\Bigl(\frac{r}{2}\Bigr)Z,\
U_{P\eta}=\sech\Bigl(\frac{r}{2}\Bigr)P\eta
    +\sech\Bigl(\frac{r}{2}\Bigr)\tanh\Bigl(\frac{r}{2}\Bigr)
    \sin^2(\varphi)Z,\
U_{PF\eta}=\sech\Bigl(\frac{r}{2}\Bigr)PF\eta.
\]
Using \eqref{expsolvprod}, \eqref{expsolv}
and \eqref{geodesic} we get
\begin{eqnarray*}
\chi_X(t)
&=& \left(\ln\sech^2\Bigl(\frac{r}{2}\Bigr),
    \Exp_{\g{n}}\left(2\tanh\Bigl(\frac{r}{2}\Bigr)\eta\right)\right)
    \cdot\Bigl(0,\Exp_{\g{n}}\left(t X\right)\Bigr)\\
&=& \left(\ln\sech^2\Bigl(\frac{r}{2}\Bigr),
    \Exp_{\g{n}}\left(2\tanh\Bigl(\frac{r}{2}\Bigr)\eta
    +t U_X\right)\right).
\end{eqnarray*}
On the other hand, we have $\chi(t):=\chi_X(t)=\exp_{g(t)\cdot
o}\bigl(s(t)(j(t)\bar{F}\eta+h(t)\eta)\bigr)$. We may write
$g(t)=(b(t),\Exp_{\g{n}}(x(t)Z+V(t))$ for certain functions
$b,x:I\to\R$ and $V:I\to\g{v}_0$. Taking into account that $g(t)$ is
an isometry, we get, using again \eqref{expsolvprod} and \eqref{geodesic},

\begin{eqnarray*}
\chi(t)\! &=& \exp_{g(t)\cdot
o}\bigl(s(t)(j(t)\bar{F}\eta+h(t)\eta)\bigr)
\ =\ g(t)\exp_{o}\bigl(s(t)(j(t)\bar{F}\eta+h(t)\eta)\bigr)\\
&=& (b(t),\Exp_{\g{n}}(x(t)Z+V(t)))
    \cdot\left(\ln\sech^2\frac{s(t)}{2},
    \Exp_{\g{n}}\left(2\tanh\frac{s(t)}{2}
    (j(t)\bar{F}\eta+h(t)\eta)\right)\right)\\
&=& \Bigl(b(t)+\ln\sech^2\frac{s(t)}{2},
    \Exp_{\g{n}}\Bigl(V(t)+2e^{b(t)/2}\tanh\frac{s(t)}{2}
    (j(t)\bar{F}\eta+h(t)\eta)\\
&&  \phantom{\Bigl(b(t)+\ln\sech^2\frac{s(t)}{2},}
    +\left\{x(t)+e^{b(t)/2}\tanh\frac{s(t)}{2}
    \langle JV(t),j(t)\bar{F}\eta+h(t)\eta\rangle
    \right\}Z\Bigr)\Bigr).
\end{eqnarray*}
As $\Exp_{\g{n}} : {\g{n}} \to N$ is a diffeomeorphism,
we easily get $j(t)\tanh(s(t)/2)=0$
for all $t$ from the previous two equations by comparing the $\bar{F}\eta$-component.
This eventually implies $j'(0)=0$, and finishes the proof for the second part of (i).
\end{proof}


\bigskip
\noindent {\sc Department of Mathematics, University College,
Cork, Ireland}\\
\noindent {\sc Email:} j.berndt@ucc.ie, jc.diazramos@ucc.ie

\begin{sidewaysfigure}
Homogeneous hypersurfaces in the complex hyperbolic space

\scriptsize
\renewcommand{\tabcolsep}{5mm} 
\renewcommand{\arraystretch}{1.7} 
\begin{tabular}{@{\vline\,}c@{\,\vline\,}c@{\,\vline\,}
c@{}c@{\,\vline\,}l@{\,\vline}}
\hline Type    &   Group acting   &   Principal curvatures
&   Multiplicities  &
    \qquad Comments\\
\hline \hline (A) &   $S(U(1,k)U(n-k))$  &
    $\frac{1}{2}\tanh\frac{r}{2}$ &   $2(n-k-1)$ &
    Tubes around totally geodesic $\CH[k]$,\\
&&  $\frac{1}{2}\coth\frac{r}{2}$ &   $2k$  &
    $0 \leq k \leq n-1$. \\
&&  $\coth r$ &   $1$   &
    Two principal curvatures if $k\in\{0,n-1\}$.\\
\hline \hline (B) &   $SO^0(1,n)$ &
    $\lambda_1=\frac{1}{2}\tanh\frac{r}{2}$ &   $n-1$ &
    Tubes around totally geodesic $\R H^n$.\\
&&  $\lambda_2=\frac{1}{2}\coth\frac{r}{2}$ &   $n-1$ &
    If $r=\ln(2+\sqrt{3})$ then $\lambda_2=\lambda_3$.\\
&&  $\lambda_3=\tanh r$ &   $1$   &\\
\hline \hline (H) &   $N$ (nilpotent part of&
    $1/2$ &   $2(n-1)$  &
    Horosphere foliation.\\
    &  Iwasawa decomposition)
    &  $1$ &   $1$ &\\
\hline \hline (S) &   $S$ (Lie algebra of $S$: &
    $\frac{3}{4}\tanh\frac{r}{2}
    +\frac{1}{2}\sqrt{1-\frac{3}{4}\tanh^2\frac{r}{2}}$
    &   $1$ &
    Solvable foliation.\\
&   $\g{s}=\g{a}+\g{w}+\g{g}_{2\alpha}$ with
    &  $\frac{3}{4}\tanh\frac{r}{2}
    -\frac{1}{2}\sqrt{1-\frac{3}{4}\tanh^2\frac{r}{2}}$
    &   $1$ &\\
&   $\g{w}$ hyperplane in $\g{g}_{\alpha}$)
    &  $\frac{1}{2}\tanh\frac{r}{2}$ &   $2n-3$    &\\
\hline \hline (W)$_{\pi/2}$ &   $N_K^0(S)S$ (Lie algebra of $S$: &
    $\lambda_1=\frac{3}{4}\tanh\frac{r}{2}
    +\frac{1}{2}\sqrt{1-\frac{3}{4}\tanh^2\frac{r}{2}}$
    &   $1$ &
    Tubes around $W^{2n-k}_{\pi/2}$, $2 \leq k \leq n-1$.\\
&   $\g{s}=\g{a}+\g{w}+\g{g}_{2\alpha}$ with
    &  $\lambda_2=\frac{3}{4}\tanh\frac{r}{2}
    -\frac{1}{2}\sqrt{1-\frac{3}{4}\tanh^2\frac{r}{2}}$
    &   $1$ &
    If $r=\ln(2+\sqrt{3})$ then $\lambda_2=\lambda_4$.\\
&   $\g{w}$ such that $\g{g}_{\alpha}\ominus\g{w}$ is real)
    &  $\lambda_3=\frac{1}{2}\tanh\frac{r}{2}$ &   $2n-k-2$  &\\
&&  $\lambda_4=\frac{1}{2}\coth\frac{r}{2}$ &   $k-1$ &\\
\hline (W)$_\varphi$ &   $N_K^0(S)S$  (Lie algebra of $S$: &
    $-\frac{1}{6} \left(\coth\frac{r}{2}
    \Bigl(u_{r,\varphi}^1+\frac{1}{u_{r,\varphi}^1}\Bigr)
    +\csch\frac{r}{2}\sech\frac{r}{2}
    +4\tanh\frac{r}{2}\right)$ &   $1$    &
    Tubes around $\W{2n-k}$, $2 \leq k \leq n-1$,\\
&   $\g{s}=\g{a}+\g{w}+\g{g}_{2\alpha}$ with
    &  $-\frac{1}{6} \left(\coth\frac{r}{2}
    \Bigl(u_{r,\varphi}^2+\frac{1}{u_{r,\varphi}^2}\Bigr)
    +\csch\frac{r}{2}\sech\frac{r}{2}
    +4\tanh\frac{r}{2}\right)$ &   $1$    &
     $k$ even. The number $u^i_{r,\varphi}$ is the $i$th cubic\\
&   $\g{w}$ such that $\g{g}_{\alpha}\ominus\g{w}$ has
    &  $-\frac{1}{6} \left(\coth\frac{r}{2}
    \Bigl(u_{r,\varphi}^3+\frac{1}{u_{r,\varphi}^3}\Bigr)
    +\csch\frac{r}{2}\sech\frac{r}{2}
    +4\tanh\frac{r}{2}\right)$ &   $1$    &
    root of
    $(\beta_{r,\varphi}+\sqrt{\beta_{r,\varphi}^2-4})/2$, where\\
&    constant K\"{a}hler angle
    &  $\frac{1}{2}\tanh\frac{r}{2}$ &   $2n-k-2$  &
    $\beta_{r,\varphi}=27\sin^2(\varphi)\tanh^2(r/2)\sech^4(r/2)-2$.\\
&   $0 < \varphi < \pi/2$)
    &  $\frac{1}{2}\coth\frac{r}{2}$ &   $k-2$ &
    Four principal curvatures if $k=2$.\\
\hline
\end{tabular}
\end{sidewaysfigure}


\begin{thebibliography}{9999}

\bibitem{Be89}
J.~Berndt: Real hypersurfaces with constant principal curvatures
in complex hyperbolic space, {\it J.\ Reine Angew.\ Math.} {\bf
395} (1989), 132--141.

\bibitem{Be98}
J.~Berndt: Homogeneous hypersurfaces in hyperbolic spaces, {\it
Math.\ Z.} {\bf 229} (1998), 589--600.

\bibitem{BB01}
J.~Berndt, M.~Br\"{u}ck: Cohomogeneity one actions on hyperbolic
spaces, {\it J.\ Reine Angew.\ Math.} {\bf 541} (2001), 209--235.

\bibitem{BCO03}
J.~Berndt, S.~Console, C.~Olmos: {\it Submanifolds and
holonomy}, Chapman \& Hall/CRC Research Notes in Mathematics {\bf
434}, Chapman \& Hall/CRC, Boca Raton, FL, 2003.

\bibitem{BD05}
J.~Berndt, J.C.~D\'{\i}az-Ramos: Real hypersurfaces with
constant principal curvatures in complex hyperbolic spaces,
to appear in {\it J.\ London Math.\ Soc.},
preprint arXiv:math.DG/0507361.

\bibitem{BD06}
J.~Berndt, J.C.~D\'{\i}az-Ramos: Real hypersurfaces with
constant principal curvatures in the complex hyperbolic plane,
to appear in {\it Proc.\ Amer.\ Math.\ Soc.},
preprint arXiv:math.DG/0605118.

\bibitem{BT05}
J.~Berndt, H.~Tamaru: Cohomogeneity one actions on noncompact
symmetric spaces of rank one, to appear in {\it Trans. Amer. Math.
Soc.}, preprint arXiv:math.DG/0505490. 

\bibitem{BTV95}
J.~Berndt, F.~Tricerri, L.~Vanhecke: {\it Generalized
Heisenberg groups and Damek--Ricci harmonic spaces}, Lecture Notes
in Mathematics {\bf 1598}, Springer--Verlag, Berlin, 1995.

\bibitem{Ca38}
\'{E}.\ Cartan: Familles de surfaces isoparam\'etriques dans les
espaces \`{a} courbure constante, {\it Ann.\ Mat.\ Pura Appl., IV.
Ser.} {\bf 17} (1938), 177--191.

\bibitem{HL71} W.-Y.~Hsiang and H.B.~Lawson~Jr.: Minimal
submanifolds of low cohomogeneity, {\it J.\ Differential Geom.}
{\bf 5} (1971), 1--38.

\bibitem{Lo98}
M.~Lohnherr:
{\it On ruled real hypersurfaces of complex space forms},
PhD Thesis, University of Cologne, 1998.

\bibitem{Ta73}
R.~Takagi: On homogeneous real hypersurfaces in a complex projective
space, {\it Osaka J.\ Math.} {\bf 10} (1973), 495--506.

\end{thebibliography}
\end{document}